\documentclass[journal,onecolumn,draftclsnofoot,10pt]{IEEEtran}

\usepackage{multirow}
\usepackage{graphicx}
\usepackage{amsmath}
\usepackage{amsthm}
\usepackage{mathtools}
\usepackage[psamsfonts]{amssymb}
\usepackage{amsxtra}
\usepackage{hyperref}
\usepackage{paralist}
\usepackage{threeparttable}
\usepackage{cleveref}
\usepackage{cite}
\usepackage{color}
\usepackage[justification=centering]{caption}

\newcommand{\argmin}{\arg\!\min}
\newtheorem{Theorem}{Theorem}[section]
\newtheorem{Lemma}[Theorem]{Lemma}
\newtheorem{proposition}[Theorem]{Proposition}

\newtheorem{corollary}[Theorem]{Corollary}
\newtheorem{Remark}[Theorem]{Remark}
\theoremstyle{definition}
\newtheorem{assump}{Assumption}

\newenvironment{myassump}[2][]
{\renewcommand\theassump{#2}\begin{assump}[#1]}
	{\end{assump}}

\newcommand{\intr}{\operatorname{int}}

\def \argmin{\mbox{argmin}}

\def \btheta{\boldsymbol{\theta}}
\def \wbtheta{\widehat{\boldsymbol{\theta}}}
\def \wx{\widehat{\mathbf{x}}}

\begin{document}
\title{Distributed Constrained Recursive Nonlinear Least-Squares Estimation: Algorithms and Asymptotics}
\author{Anit~Kumar~Sahu,~\IEEEmembership{Student Member,~IEEE}, Soummya~Kar,~\IEEEmembership{Member,~IEEE},\\
Jos\'e M.~F.~Moura,~\IEEEmembership{Fellow,~IEEE} and H.~Vincent~Poor,~\IEEEmembership{Fellow,~IEEE}
\thanks{
Digital Object Identifier: 10.1109/TSIPN.2016.2618318

Copyright (c) 2016 IEEE. Personal use of this material is permitted. However, permission to use this material for any other purposes must be obtained from the IEEE by sending a request to pubs-permissions@ieee.org.

A. K. Sahu, S. Kar and J. M. F. Moura are with the Department of Electrical and Computer Engineering, Carnegie Mellon University, Pittsburgh, PA 15213, USA (email:anits@andrew.cmu.edu, soummyak@andrew.cmu.edu, moura@ece.cmu.edu). H. V. Poor is with the Department of Electrical Engineering, Princeton University, Princeton, NJ 08544 USA (e-mail: poor@princeton.edu). 

The work of A. K. Sahu and S. Kar was supported in part by NSF under grants CCF-1513936 and ECCS-1408222. The work of J. M. F. Moura was supported in part by NSF under grants CCF-1011903 and CCF-1513936. The work of H. V. Poor was supported in part by NSF under grant CCF-1420575 and by ARO under MURI grant W911NF-11-1-0036.
}}

\maketitle

\begin{abstract}
This paper focuses on recursive nonlinear least squares parameter estimation in multi-agent networks, where the individual agents observe sequentially over time an independent and identically distributed (i.i.d.) time-series consisting of a nonlinear function of the true but unknown parameter corrupted by noise. A distributed recursive estimator of the \emph{consensus}+\emph{innovations} type, namely $\mathcal{CIWNLS}$, is proposed, in which the agents update their parameter estimates at each observation sampling epoch in a collaborative way by simultaneously processing the latest locally sensed information~(\emph{innovations}) and the parameter estimates from other agents~(\emph{consensus}) in the local neighborhood conforming to a pre-specified inter-agent communication topology. Under rather weak conditions on the connectivity of the inter-agent communication and a \emph{global observability} criterion, it is shown that, at every network agent, $\mathcal{CIWNLS}$ leads to consistent parameter estimates. Furthermore, under standard smoothness assumptions on the local observation functions, the distributed estimator is shown to yield order-optimal convergence rates, i.e., as far as the order of pathwise convergence is concerned, the local parameter estimates at each agent are as good as the optimal centralized nonlinear least squares estimator that requires access to all the observations across all the agents at all times. To benchmark the performance of the $\mathcal{CIWNLS}$ estimator with that of the centralized nonlinear least squares estimator, the asymptotic normality of the estimate sequence is established, and the asymptotic covariance of the distributed estimator is evaluated. Finally, simulation results are presented that illustrate and verify the analytical findings.
\end{abstract}

\begin{IEEEkeywords}
Multi-agent networks, distributed estimation, nonlinear least squares, distributed information processing, \emph{consensus}+\emph{innovations}, distributed stochastic aproximation, distributed inference, networked inference
\end{IEEEkeywords}

\section{Introduction}
\label{introduction}

The paper focuses on distributed nonlinear least squares estimation in distributed information settings. Each agent in the network senses sequentially over time independent and identically distributed (i.i.d) time-series that are (nonlinear) functions of the underlying vector parameter of interest corrupted by noise. To be specific, we are interested in the design of recursive estimation algorithms to estimate a vector parameter of interest that are consistent and order-optimal in the sense of pathwise convergence rate and such that their asymptotic error covariances are comparable with that of the centralized weighted nonlinear least squares estimator\footnote{A centralized estimator has access to all agent data at all times and has sufficient computing ability to implement the classical weighted nonlinear least squares estimator~\cite{Jeinrich,wu1981asymptotic} at all times.}. The estimation algorithms we design are recursive -- they process the agents' observations at all times as and when they are sensed, rather than batch processing. This contrasts with centralized setups, where a fusion center has access to all the observations across different agents at all times, i.e., the inter-agent communication topology is all-to-all or all-to-one. Centralized estimators are burdened by high communication overheads, synchronization issues, and high energy requirements. Moreover, there is the requirement of global model information, i.e., the fusion center requiring information about the local models of all agents. All these make centralized estimation algorithms difficult to implement in multi-agent distributed setups of the type considered in this paper, motivating us to revisit the problem of distributed \emph{sequential} parameter estimation. To accommodate energy constraints in many practical networked and wireless settings, the inter-agent collaboration is limited to a pre-assigned possibly sparse communication graph. Moreover, due to limited computation and storage capabilities of individual agents in a typical multi-agent networked setting, we restrict to scenarios where individual agents are only aware of their local model information; hence, we allow for heterogeneity among agents, with different agents possibly having different local sensing models and noise statistics.
This paper proposes a distributed recursive algorithm, namely, the $\mathcal{CIWNLS}$~($\emph{Consensus}+\emph{innovations}$ Weighted Nonlinear Least Squares), which is of the $\emph{consensus}+\emph{innovations}$ form~\cite{KarMouraRamanan-Est-2008}. We specifically focus on a setting in which the agents make i.i.d observations sequentially over time, only possess local model information, and update their parameter estimates by simultaneous assimilation of the information obtained from their neighboring agents~(\emph{consensus}) and current locally sensed information~(\emph{innovation}). This justifies the name $\mathcal{CIWNLS}$, which is a distributed weighted nonlinear least squares~(WNLS) type algorithm of the $\emph{consensus}+\emph{innovations}$ form. To replicate practical sensing environments accurately, we model the underlying vector parameter as a static parameter, that takes values in a parameter set $\Theta\subseteq\mathbb{R}^{M}$~(possibly a strict subset of $\mathbb{R}^{M}$).  The dimension $M$ is possibly large, but the observation of any agent $n$ is $M_{n}$ dimensional with typically $M_{n}\ll M$ in most applications; this renders the parameter locally unobservable at each agent. The key assumptions concerning the sensing functions in this paper are required to hold only on the parameter set $\Theta$ and not on the entire space\footnote{By taking the parameter set $\Theta=\mathbb{R}^{M}$, the unconstrained parameter estimation problem can be addressed, and thus the setup in this paper enables a richer class of formulations.} $\mathbb{R}^{M}$. The distributed sequential estimation approach of the $\emph{consensus}+\emph{innovations}$ form that we present accomplishes the following:

\noindent {\bf Consistency under global observability}: We assume \emph{global observability}\footnote{Global observability corresponds to the centralized setting, where an estimator has access to the observations of all sensors at all times. The assumption of global observability does not mean that each sensor is observable; rather, if there was a centralized estimator with simultaneous access to all the sensor measurements, this centralized estimator would be able to reasonably estimate the underlying parameter. A more precise definition is provided later in Assumption \ref{m:1}.} and certain \emph{monotonicity} properties of the multi-agent sensing model, as well as the connectedness of the inter-agent communication graph. We show that our recursive distributed estimator generates parameter estimate sequences that are strongly consistent~(see, \cite{wasserman2013all,le2012asymptotics} for a detailed treatment on consistency) at each agent. Global observability is a minimal requirement for consistency; in fact, it is necessary for consistency of centralized estimators as well.

\noindent {\bf Optimal pathwise convergence rate}\footnote{By \emph{optimal pathwise convergence rate}, we mean the pathwise convergence rate of the centralized estimator to the true underlying parameter with noisy observations. }: We show that the proposed distributed estimation algorithm $\mathcal{CIWNLS}$ yields order-optimal pathwise convergence rate under certain smoothness conditions on the sensing model. These conditions are standard in the recursive estimation literature and we require them to hold only on the parameter set $\Theta$. Even though recursive, our distributed estimation approach guarantees that the parameter estimates are feasible at all times, i.e., they belong to the parameter set $\Theta$. Further, the parameter estimates at each local agent $n$ are as good as the optimal centralized estimator as far as pathwise convergence rate is concerned. The key point to note here is that, for the above order optimality to hold we need to only assume that the inter-agent communication graph is connected irrespective of how sparse the link realizations are.

\noindent {\bf Asymptotic Normality}: Under standard smoothness conditions on the sensing model, the proposed distributed estimation algorithm $\mathcal{CIWNLS}$ is shown to yield asymptotically normal\footnote{An estimate sequence is asymptotically normal if its $\sqrt{t}$ scaled error process, i.e., the difference between the sequence and the true parameter converges in distribution to a normal random variable, where $t$ refers to (discrete) time or equivalently the number of sampling epochs.} parameter estimate sequences~(see, \cite{wasserman2013all,le2012asymptotics} for a detailed treatment on asymptotic normality). Distributed estimation does pay a price. The asymptotic covariance of the proposed distributed estimator is not as efficient as that of the centralized estimator; nonetheless, it shows the benefits of inter-agent collaboration. In absence of inter-agent collaboration, the parameter of interest most likely is unobservable at each individual agent, and hence non-collaborative or purely decentralized procedures will lead to divergence under the usual asymptotic normality scaling at the individual network agents.

\noindent{\bf Related Work}: Distributed inference approaches addressing problems related to distributed estimation, parallel computing, and optimization in multi-agent environments through interacting stochastic gradient and stochastic approximation algorithms have been developed extensively in the literature~--see, for example, early work~\cite{tsitsiklisphd84,tsitsiklisbertsekasathans86,Bertsekas-survey,Kushner-dist}. Existing distributed estimation schemes in the literature can be broadly divided into three classes. The first class includes architectures that are characterized by the presence of a fusion center~(see, for example \cite{aysal2008constrained,luo2005universal}) that receives the estimates or local measurements or their quantized versions from the network agents and performs estimation. The second class involves single snapshot data collection~(see, for example \cite{das2006distributed,schizas2008consensus}) followed by distributed consensus or optimization protocols to fuse the initial estimates. In contrast to these classes, the third class involves agents making observations sequentially over time and where inter-agent communication, limited to arbitrary pre-assigned possibly sparse topologies, occurs at the same rate as sensing~(see, for example \cite{braca2010asymptotic,braca2008enforcing,KarMouraRamanan-Est-2008,Sayed-LMS}). Two representative schemes from the third class are \emph{consensus}+\emph{innovations} type ~\cite{KarMouraRamanan-Est-2008,Kar-Moura-LinEst-Asilomar-2008} and diffusion type algorithms~\cite{Sayed-LMS,chen2012diffusion,matta2014diffusion,cattivelli2011distributed,cattivelli2010diffusion}. Broadly speaking, these algorithms simultaneously assimilate a single round of neighborhood information, consensus like in~\cite{Bertsekas-survey,olfatisaberfaxmurray07,Dimakis-Gossip-SPM-2011,jadbabailinmorse03}, with the locally sensed latest information, the local innovation; see for example
\emph{consensus}+\emph{innovation} approaches for nonlinear distributed estimation~\cite{KarMouraRamanan-Est-2008,kar2013consistent} and detection~\cite{bajovicjakoveticxaviresinopolimoura-11,jakovetic2012distributed,sahu2015glrttit}. A key difference between the diffusion algorithms discussed above and the \emph{consensus}+\emph{innovations} algorithms presented in this paper is the nature of the innovation gains (the new information fusion weights). In the diffusion framework, the innovation gains are taken to be constant, whereas, in the \emph{consensus}+\emph{innovations} schemes these are made to decay over time in a controlled fashion. The constant innovation gains in the diffusion approaches facilitate adaptation in dynamic parameter environments, but at the same time lead to non-zero residual estimation error at the agents~(see, for example, \cite{Sayed-LMS}), whereas the time-varying innovation weights in the \emph{consensus}+\emph{innovations} approach ensure consistent parameter estimates at the agents. We also note that \cite{mateos2009distributed,mateos2012distributed} consider the problem of distributed recursive least squares using constant step-size recursive distributed algorithms for the estimate update. In order to ensure that their proposed algorithm is adaptive, the weights~(the innovation gains) in their algorithm are static, but the adaptivity comes at a loss in terms of accuracy, i.e., the algorithms result in asymptotic residual non-zero mean square error. The observation model considered in~\cite{mateos2009distributed,mateos2012distributed} is linear. In comparison, we consider a general non-linear observation model and propose an algorithm of the $\emph{consensus}+\emph{innovations}$ form, in which the innovation gains are made to decay in a controlled manner so as to ensure consistency and almost sure convergence of the sequence of parameter estimates. Other approaches for distributed inference in multi-agent networks have been considered, see for example algorithms for network inference and optimization, networked LMS and variants~\cite{Stankovic-parameter,Giannakis-LMS,Nedic-parameter,weng2014efficient,olfati2007distributed,hlinka2013distributed}. \\
\noindent More recently, in~\cite{kar2011convergence,Kar-AdaptiveDistEst-SICON-2012,kar2014asymptotically} asymptotically efficient (in the sense of optimal asymptotic covariance and optimal decay of estimation errors) \emph{consensus}+\emph{innovations} estimation procedures were presented for a wide class of distributed parameter estimation scenarios including nonlinear models and settings with imperfect statistical model information. In the context of the current paper, as the parameter belongs to a constrained set, we consider in addition a local projection to ensure that the parameter estimate is feasible at all times. Distributed iterative algorithms that include a (local) projection step have been proposed~(see, for example,~\cite{ram2010distributed}) to ensure convergence in the context of distributed optimization.\\
The \emph{consensus}+\emph{innovations} procedures generalize stochastic approximation~(see \cite{robbins1951stochastic} for an early work) to distributed multi-agent networked settings. To achieve optimality~(for instance, optimal decay of errors in a parameter estimation setup), these algorithms are designed to have a mixed time-scale flavor. By mixed time-scale, we mean that the inter-agent communication~(consensus) occurs at the same rate as that of observation sampling or incorporation of the latest sensed information~(innovation), however, the consensus and innovation terms in the iterative update are weighed by two different weight sequences, which decay to zero in a controlled manner at different rates.~(It is to be noted that this mixed time-scale is different from stochastic approximation algorithms with coupling~(see \cite{Borkar-stochapp}), where a quickly switching parameter influences the relatively slower dynamics of another state, leading to averaged dynamics.)  The \emph{consensus}+\emph{innovations} class of algorithms not only converges to the true value of the underlying parameter, but also yields asymptotic efficiency, see, for example~\cite{kar2011convergence}. We also note that in \cite{Gelfand-Mitter}, in the context of optimization, methods pertaining to mixed time-scale stochastic approximation algorithms are developed, albeit in the centralized context. The corresponding innovation term in the algorithm proposed in \cite{Gelfand-Mitter} is a martingale difference term. However, in \emph{consensus}+\emph{innovations} algorithms, see \cite{kar2011convergence,KarMouraRamanan-Est-2008}, and this work the innovation term is not a martingale difference sequence and hence it is of particular interest to characterize the rate of convergence of the innovation sequence to a martingale difference sequence, so as to establish convergence and consistency of the parameter estimate sequence. \\
\noindent We contrast the current algorithm $\mathcal{CIWNLS}$ with the distributed estimation algorithm of the \emph{consensus}+\emph{innovations} form developed in \cite{kar2014asymptotically} for a very general nonlinear setup. In \cite{kar2014asymptotically}, strong consistency of the parameter estimate sequence is established, and it is shown that the proposed algorithm is asymptotically efficient, i.e., its asymptotic covariance is the same as that of the optimal centralized estimator. However, in \cite{kar2014asymptotically}, the smoothness assumptions on the sensing functions need to hold on the entire parameter space, i.e., $\mathbb{R}^{M}$. In contrast, we consider here a setup where the parameter belongs to a constrained set and the smoothness conditions on the sensing functions need to hold only on the constrained parameter set; this allows the algorithm proposed in this paper, namely $\mathcal{CIWNLS}$, to be applicable to other types of application scenarios. Moreover, in \cite{kar2014asymptotically} the problem setup needs more detailed knowledge of the statistics of the noise processes involved, as it aims to obtain asymptotically efficient (in that the agent estimates are asymptotically normal with covariance equal to the inverse of the associated Fisher information rate) estimates for general statistical exponential families. In particular, to achieve asymptotic efficiency, \cite{kar2014asymptotically} develops a \emph{consensus}+\emph{innovations} type distributed recursive variant of the maximum likelihood estimator (MLE) that requires knowledge of the detailed observation statistics. In contrast, in this paper, our setup only needs knowledge of the noise covariances and the sensing functions. Technically speaking, for additive noisy observation models, the weighted nonlinear squares estimation, the distributed version of which is proposed in this paper, applies to fairly generic estimation scenarios, i.e., where observation noise statistics are unknown. In a previous work \cite{kar2013consistent}, we considered a similar setup as used in this paper and proposed a \emph{consensus}+\emph{innovations} type distributed estimation algorithm, where the parameter set is a compact convex subset of the $M$-dimensional Euclidean space and established the consistency and order-optimality of the parameter estimate sequence.  In this paper, we not only establish the consistency and order-optimality of the parameter estimate sequence, but also the asymptotic normality of the parameter estimate sequence.\\
\noindent Distributed inference algorithms have been applied to various other problems in the networked setting, such as, distributed optimization and distributed detection. Distributed optimization in networked settings has been considered, see for example~\cite{nedic2015decentralized,koppel2015saddle,tsianos2012consensus,sundhar2010distributed}~that established that distributed gradient or distributed subgradient based algorithms converge to the optimum solutions for a general class of optimization problems. \emph{Consensus}+\emph{innovations} algorithms~(see, \cite{Kar-Moura-LinEst-Asilomar-2008} for example) have been proposed to address problems pertaining to networked inference and networked optimization~(see, \cite{jakovetic2014convergence,kar2014asymptotically} for example). Distributed detection has been very popular of late, where distributed inference algorithms of the third class, as discussed above, which assimilate information from neighbors and current locally sensed observations at the same rate, have been extensively used. References \cite{bajovicjakoveticxaviresinopolimoura-11,jakovetic2012distributed,sahu2015glrttit} use large deviations theory to find decay rates of error probabilities for distributed detectors of the \emph{consensus}+\emph{innovations} type; see also subsequent work~\cite{matta2015exact,matta2016distributed} that study similar problems for the class of diffusion algorithms.

\noindent\textbf{Paper Organization:}
The rest of the paper is organized as follows. The notation to be used throughout the paper is presented in Section \ref{subsec:not}, where spectral graph theory is also reviewed. The multi-agent sensing model is described in Section~\ref{sys-model}, where we also review some classical concepts on estimation theory. Section~\ref{alg_DQ} presents the proposed distributed parameter estimation algorithm $\mathcal{CIWNLS}$. The main results of this paper concerning the consistency and the asymptotic normality of the parameter estimate sequence are provided in Section~\ref{sec:main_res}. Section \ref{sec:sim} presents the simulation results. The proof of the main results of this paper are provided in Section \ref{sec:proof_main_res}. Finally, Section~\ref{conclusion} concludes the paper and discusses future research avenues.\\
\subsection{Notation}
\label{subsec:not}
We denote by~$\mathbb{R}$ the set of reals, by~$\mathbb{R}_{+}$ the set of non-negative reals, and by~$\mathbb{R}^{k}$ the $k$-dimensional Euclidean space.
\noindent The set of $k\times k$ real matrices is denoted by $\mathbb{R}^{k\times k}$. The set of integers is  $\mathbb{Z}$, whereas, $\mathbb{Z}_{+}$ is the subset of non-negative integers. Vectors and matrices are in bold faces; $\mathbf{A}_{ij}$ or $[\mathbf{A}]_{ij}$ the $(i,j)$-th entry of a matrix $\mathbf{A}$; $\mathbf{a}_{i}$ or $[\mathbf{a}]_{i}$ the $i$-th entry of a vector $\mathbf{a}$. The symbols $\mathbf{I}$ and $\mathbf{0}$ are the $k\times k$ identity matrix and the $k\times k$ zero matrix, respectively, the dimensions being clear from the context. The vector $\mathbf{e_{i}}$ is the $i$-th column of $\mathbf{I}$. The symbol $\top$ stands for matrix transpose. The determinant and trace of a matrix are $\det(.)$ and $\operatorname{tr}(.)$, respectively. The $k\times k$ matrix $\mathbf{J}=\frac{1}{k}\mathbf{1}\mathbf{1^{\top}}$, where $\mathbf{1}$ is the $k\times 1$ vector of ones. The operator $\otimes$ denotes the Kronecker product. The operator $|| . ||$ applied to a vector is the standard Euclidean $\mathcal{L}_{2}$ norm, while when applied to matrices stands for the induced $\mathcal{L}_{2}$ norm, which is equivalent to the spectral radius for symmetric matrices. The cardinality of a set $\mathcal{S}$ is $\left|\mathcal{S}\right|$.\\
Throughout the paper, we assume that all random objects are defined on a common
measurable space $(\Omega, \mathcal{F})$ equipped with a filtration (sequence of increasing sub-$\sigma$-algebras of $\mathcal{F}$) $\{\mathcal{F}_{t}\}$~(see~\cite{Jacod-Shiryaev}). The true (but unknown) value of the parameter is denoted by $\btheta$. For the true but unknown parameter value $\btheta$, probability and expectation on the measurable space $(\Omega,\mathcal{F})$ are written as $\mathbb{P}_{\btheta} \left[\cdot\right]$ and $\mathbb{E}_{\btheta}\left[\cdot\right]$, respectively. A stochastic process $\{\mathbf{z}_{t}\}$ is said to be adapted to a filtration $\{\mathcal{G}_{t}\}$ or $\{\mathcal{G}_{t}\}$-adapted if $\mathbf{z}_{t}$ is measurable with respect to (w.r.t.) the $\sigma$-algebra $\mathcal{G}_{t}$ for all $t$. All inequalities involving random variables are to be interpreted almost surely~(a.s.). \\
\noindent For deterministic $\mathbb{R}_{+}$-valued sequences $\{a_{t}\}$ and $\{b_{t}\}$, the notation $a_{t}=O(b_{t})$ implies the existence of a constant $c>0$ such that $a_{t}\leq cb_{t}$ for all $t$ sufficiently large; the notation $a_{t}=o(b_{t})$ denotes $a_{t}/b_{t}\rightarrow 0$ as $t\rightarrow\infty$. The order notations $O(\cdot)$ and $o(\cdot)$ will be used in the context of stochastic processes as well in which case they are to be interpreted almost surely or pathwise.

\noindent\textbf{Spectral graph theory}: The inter-agent communication network is an \emph{undirected} simple connected graph $G=(V,E)$, with $V=\left[1\cdots N\right]$ and~$E$ denoting the set of agents (nodes) and communication links, see~\cite{SensNets:Bollobas98}. The neighborhood of node~$n$ is
\begin{equation}
\label{def:omega} \Omega_{n}=\left\{l\in V\,|\,(n,l)\in
E\right\}.
\end{equation}
The node~$n$ has degree $d_{n}=|\Omega_{n}|$. The structure of the graph is described by the  $N\times N$ adjacency matrix, $\mathbf{A}=\mathbf{A}^\top=\left[\mathbf{A}_{nl}\right]$, $\mathbf{A}_{nl}=1$, if $(n,l)\in E$, $\mathbf{A}_{nl}=0$, otherwise. Let $\mathbf{D}=\mbox{diag}\left(d_{1}\cdots d_{N}\right)$. The graph Laplacian $\mathbf{L}=\mathbf{D}-\mathbf{A}$ is positive definite, with eigenvalues ordered as $0=\lambda_{1}(\mathbf{L}) \leq\lambda_{2}(\mathbf{L}) \leq\cdots \leq \lambda_{N}(\mathbf{L})$. The eigenvector of~$\mathbf{L}$ corresponding to $\lambda_{1}(\mathbf{L})$ is $(1/\sqrt{N})\mathbf{1}_{N}$. The multiplicity of its zero eigenvalue equals the number of connected components of the network; for a connected graph, $\lambda_{2}(\mathbf{L})>0$. This second eigenvalue is the algebraic connectivity or the Fiedler value of the network~(see~\cite{chung1997spectral}~for instance).

\section{Sensing Model and Preliminaries}
\label{sys-model}
\noindent Let $\btheta\in\Theta$~(to be specified shortly) be an $M$-dimensional (vector) parameter that is to be estimated by a network of $N$ agents.
\noindent We specifically consider a discrete time system. Each agent $n$ at time $t$ makes a noisy observation $\mathbf{y}_{n}(t)$ that is a noisy function (nonlinear) of the parameter. \noindent Formally, the observation model for the $n$-th agent is given by
\begin{equation}
\label{eq:obsmod}
\mathbf{y}_{n}(t)=\mathbf{f}_{n}(\btheta)+\mathbf{\zeta}_{n}(t),
\end{equation}
\noindent where $\mathbf{f}_{n}(\cdot)$ is, in general, a non-linear function, $\{\mathbf{y}_{n}(t)\}$ is a $\mathbb{R}^{M_{n}}$-valued observation sequence for the $n$-th agent and for each $n$, $\left\{\mathbf{\zeta}_{n}(t)\right\}$ is a zero-mean temporally independent and identically distributed~(i.i.d.)~ noise sequence with nonsingular covariance matrix $\mathbf{R}_{n}$, such that, $\mathbf{\zeta}_{n}(t)$ is $\mathcal{F}_{t+1}$-adapted and independent of $\mathcal{F}_{t}$.
\noindent In typical application scenarios, the observation at each agent is low-dimensional, i.e., $M_{n}\ll M$, and usually a function of only a subset
of the $M$ components of $\btheta$, i.e., agent $n$ observes a function of $K_{n}$ components of $\btheta$ with $K_{n}\ll M$, which most likely entails that the parameter of interest $\btheta$ is locally unobservable at the individual agents. Hence, to achieve a reasonable estimate of the parameter $\btheta$, it is necessary for the agents to collaborate through inter-agent message passing schemes.

\noindent Since, the sources of randomness in our formulation are the observations $\mathbf{y}_{n}(t)$'s by the agents, the filtration $\{\mathcal{F}_{t}\}$ (introduced in Section \ref{subsec:not})~may be taken to be the natural filtration generated by the random observations,~i.e.,
\begin{align}
 \label{eq:filt_1}
 \mathcal{F}_{t}=\mathbf{\sigma}\left(\left\{\left\{\mathbf{y}_{n}(s)\right\}_{n=1}^{N}\right\}_{s=0}^{t-1}\right),
\end{align}
\noindent which is the $\sigma$-algebra induced by the observation processes.

\noindent To motivate our distributed estimation approach (presented in Section~\ref{alg_DQ}) and benchmark its performance with respect to the optimal centralized estimator, we now review some concepts from centralized estimation theory.\\
\noindent{\bf Centralized weighted nonlinear least-squares (WNLS) estimation}. Consider a network of agents where a hypothetical fusion center has access to the observations made by all the agents at all times and then conducts the estimation scheme. In such scenarios, one of the most widely used estimation approaches is the weighted nonlinear least squares~(WNLS)~(see, for example,~\cite{Jeinrich}). The WNLS is applicable to fairly generic estimation scenarios, for instance, even when the observation noise statistics are unknown, which precludes other classical estimation approaches such as the maximum likelihood estimation. We discuss below some useful theoretical properties of WNLS like, in case the observation noise is Gaussian, it coincides with the (asymptotically) efficient maximum likelihood estimator. To formalize, for each $t$, define the cost function
\begin{equation}
\label{def_costQ}
\mathcal{Q}_{t}\left(\mathbf{z}\right)=\sum_{s=0}^{t}\sum_{n=1}^{N}\left(\mathbf{y}_{n}(s)-\mathbf{f}_{n}(\mathbf{z})\right)^{\top}\mathbf{R}_{n}^{-1}\left(\mathbf{y}_{n}(s)-\mathbf{f}_{n}(\mathbf{z})\right),
\end{equation}
where $\mathbf{R}_{n}$ denotes the positive definite covariance of the measurement noise $\mathbf{\zeta}_{n}(t)$. The WNLS estimate $\wbtheta_{t}$ of $\btheta$ at each time $t$ is obtained by minimizing the cost functional $\mathcal{Q}_{t}(\cdot)$,
\begin{equation}
\label{WNLS}
\wbtheta_{t}\in\argmin_{\mathbf{z}\in\Theta}\mathcal{Q}_{t}(\mathbf{z}).
\end{equation}
Under rather weak assumptions on the sensing model (stated below), the existence and asymptotic behavior of WNLS estimates have been analyzed in the literature.\\
\begin{myassump}{M1}
	\label{m:2}
	\emph{The set $\Theta$ is a closed convex subset of $\mathbb{R}^{M}$ with non-empty interior $\intr(\Theta)$ and the true \textup{(}but unknown\textup{)} parameter $\btheta\in\intr(\Theta)$}.\\
\end{myassump}
\begin{myassump}{M2}
\label{m:1}
\emph{The sensing model is globally observable, i.e., any pair $\btheta, \acute{\btheta}$ of possible parameter instances in $\Theta$ satisfies
	\begin{equation}
	\label{WNLS1}
	\sum_{n=1}^{N}\left\|\mathbf{f}_{n}(\btheta)-\mathbf{f}_{n}(\acute{\btheta})\right\|^{2}=0
	\end{equation}
	if and only if $\btheta=\acute{\btheta}$.}\\
\end{myassump}

\begin{myassump}{M3}
	\label{m:3}
	\emph{The sensing function $\mathbf{f}_{n}(.)$ for each $n$ is continuously differentiable in the interior $\intr(\Theta)$ of the set $\Theta$. For each $\btheta$ in the set $\Theta$, the matrix $\mathbf{\Gamma}_{\btheta}$ that is given by
		\begin{align}
		\label{eq:m:3_1}
		\mathbf{\Gamma}_{\btheta}=\frac{1}{N}\sum_{n=1}^{N}\nabla\mathbf{f}_{n}\left(\btheta\right)\mathbf{R}_{n}^{-1}\nabla\mathbf{f}_{n}^{\top}\left(\btheta\right),
		\end{align}
		where $\nabla \mathbf{f}$ denotes the gradient of $\mathbf{f}(\cdot)$, is invertible.}\\
\end{myassump}
\noindent Smoothness conditions on the sensing functions, such as the one imposed by assumption \ref{m:3} is common in the literature addressing statistical inference algorithms in non-linear settings.
\noindent Note that the matrix $\mathbf{\Gamma}_{\btheta}$ is well defined at the true value of the parameter $\btheta$ as $\btheta\in\intr(\Theta)$ and the continuous differentiability of the sensing functions hold for all $\btheta\in\intr(\Theta)$.

\begin{myassump}{M4}
\label{m:4}
\emph{There exists $\epsilon_{1}>0$, such that, for all $n$, $\mathbb{E}_{\btheta}\left[\left\|\zeta_{n}(t)\right\|^{2+\epsilon_{1}}\right]<\infty$.}
\end{myassump}

\noindent The following classical result characterizes the asymptotic properties of (centralized) WNLS estimators.
\begin{proposition}\textup{(\hspace{-0.2pt}\cite{Jeinrich})}
\label{prop:WNLS}
Let the parameter set $\Theta$ be compact and the sensing function $f_{n}(\cdot)$ be continuous on $\Theta$ for each $n$. Then, a WNLS estimator of $\btheta$ exists, i.e., there exists an $\{\mathcal{F}_{t}\}$-adapted process $\{\wbtheta_{t}\}$ such that
\begin{equation}
\label{WNLS2}
\wbtheta_{t}\in\textrm{argmin}_{\mathbf{z}\in\Theta}\mathcal{Q}_{t}(\mathbf{z}),~\forall t.
\end{equation}
Moreover, if the model is globally observable, i.e., Assumption \ref{m:1} holds, the WNLS estimate sequence $\{\wbtheta_{t}\}$ is consistent, i.e.,
\begin{equation}
\label{WNLS3}
\mathbb{P}_{\btheta}\left(\lim_{t\rightarrow\infty}\wbtheta_{t}=\btheta\right)=1.
\end{equation}
Additionally, if Assumption \ref{m:3} holds, the parameter estimate sequence is asymptotically normal, i.e.,

\begin{align}
\label{eq:asymp_var_1}
\sqrt{t+1}\left(\wbtheta_{t}-\btheta\right)\overset{D}{\Longrightarrow}\mathcal{N}\left(0, \mathbf{\Sigma}_{c}\right),
\end{align}
where
\begin{align}
\label{eq:asymp_var_11}
\mathbf{\Sigma}_{c}=\left(N\mathbf{\Gamma}_{\btheta}\right)^{-1},
\end{align}
$\mathbf{\Gamma}_{\btheta}$ is as given by \eqref{eq:m:3_1} and $\overset{\mathcal{D}}{\Longrightarrow}$ refers to convergence in distribution (weak convergence).
\end{proposition}

\noindent The WNLS estimator, apart from needing a fusion center that has access to the observations across all agents at all times, also incorporates a batch data processing as implemented in \eqref{WNLS}. To mitigate the enormous communication overhead incurred in \eqref{WNLS}, much work in the literature has focused on the development of sequential albeit centralized estimators that process the observations $\mathbf{y}(t)$ across agents in a recursive manner. Under additional smoothness assumptions on the local observation functions $\mathbf{f}_{n}(\cdot)$'s, recursive centralized estimators of the stochastic approximation type have been developed by several authors, see, for example,~\cite{Sakrison1965efficient,Hasminskii1974estimation,Pfanzagl1974efficient,Stone1975sequentialestimation,Fabian1978efficient}. Such centralized estimators require a fusion center that process the observed data in batch mode or recursive form. The fusion center receives the entire set of agents' observations, $\{\mathbf{y}_{n}(t)\}$,~$n=1, 2, 3,\cdots, N,$ at all times $t$. Moreover, both in the batch and the recursive processing form, the fusion center needs global model information in the form of the local observation functions $\mathbf{f}_{n}(\cdot)$'s and the observation noise statistics, i.e., the noise covariances $\mathbf{R}_{n}$'s across all agents. In contrast, this paper develops collaborative distributed estimators of $\btheta$ at each agent $n$ of the network, where each agent $n$ has access to its local sensed data $\mathbf{y}_{n}(t)$ only and local model information, i.e., its own local sensing function $\mathbf{f}_{n}(\cdot)$ and noise covariance $\mathbf{R}_{n}$. To mitigate the communication overhead, we present distributed message passing schemes in which agents, instead of forwarding raw data to a fusion center, participate in a collaborative iterative process to estimate the underlying parameter $\btheta$. The agents also maintain a copy of their local parameter estimate that is updated by simultaneously processing local parameter estimates from their neighbors and the latest sensed information. To obtain a good parameter estimate with such localized communication, we propose a distributed estimator that incorporates neighborhood information mixing and local data processing simultaneously (at the same rate). Such estimators are referred to as $\emph{consensus}+\emph{innovations}$ estimators, see  \cite{KarMouraRamanan-Est-2008}, for example.

\section{A Distributed Estimator : $\mathcal{CIWNLS}$}
\label{alg_DQ} \noindent We state formally assumptions pertaining to the inter-agent communication and additional smoothness conditions on the sensing functions required in the distributed setting.\\

\begin{myassump}{M5}
\label{m:5}\emph{The inter-agent communication graph is connected, i.e., $\lambda_{2}(\mathbf{L}) > 0$, where $\mathbf{L}$ denotes the associated graph Laplacian matrix.}
\end{myassump}

\begin{myassump}{M6}
	\label{m:6}\emph{For each $n$, the sensing function $\mathbf{f}_{n}(\cdot)$ is Lipschitz continuous on $\Theta$, i.e., for each agent $n$, there exists a constant $k_{n} >0$ such that
	\begin{align}
	\label{eq:as2_1}
	\left\|\mathbf{f}_{n}\left(\btheta\right)-\mathbf{f}_{n}\left(\btheta^{*}\right)\right\| \le k_{n}\left\|\btheta-\btheta^{*}\right\|,
	\end{align}
	for all $\btheta, \btheta^{*} \in \Theta$.}
\end{myassump}

\noindent{\bf Distributed algorithm}. In the proposed implementation, each agent $n$ updates at each time $t$ its estimate sequence $\{\mathbf{x}_{n}(t)\}$ and an auxiliary sequence $\{\wx_{n}(t)\}$ using a two-step collaborative procedure; specifically, 1)~$\wx_{n}(t)$ is updated by a \emph{consensus}+\emph{innovations} rule and, subsequently, 2)~a local projection to the feasible parameter set $\Theta$ updates $\mathbf{x}_{n}(t)$. Formally, the overall update rule at an agent $n$ corresponds to
\begin{align}
\label{eq:dist1}
&\wx_{n}(t+1)=\mathbf{x}_{n}(t)-\underbrace{\beta_{t}\sum_{l\in\Omega_{n}}\left(\mathbf{x}_{n}(t)-\mathbf{x}_{l}(t)\right)}_{\mathrm{neighborhood~ consensus}}\nonumber\\
&-\underbrace{\alpha_{t}\left(\nabla \mathbf{f}_{n}(\mathbf{x}_{n}(t))\right)\mathbf{R}_{n}^{-1}\left(\mathbf{f}_{n}(\mathbf{x}_{n}(t))-\mathbf{y}_{n}(t)\right)}_{\mathrm{local~ innovation}}
\end{align}
and
\begin{equation}
\label{eq:dist2}
\mathbf{x}_{n}(t+1)=\mathcal{P}_{\Theta}[\wx_{n}(t+1)],
\end{equation}
where : $\Omega_{n}$ is the communication neighborhood of agent $n$~(determined by the Laplacian $\mathbf{L}$); $\nabla f_{n}(\cdot)$ is the gradient of $\mathbf{f}_{n}$, which is a matrix of dimension $\mathbf{M}\times\mathbf{M}_{n}$, with the $(i,j)$-th entry given by $\frac{\partial \left[\mathbf{f}_{n}\left(\mathbf{x}_{n}(t)\right)\right]_{j}}{\partial \left[\mathbf{x}_{n}(t)\right]_{i}}$; $\mathcal{P}_{\Theta}[\cdot]$ the projection operator corresponding to projecting\footnote{The projection on $\Theta$ is unique under assumption~\ref{m:2}.} on $\Theta$; and $\{\beta_{t}\}$ and $\{\alpha_{t}\}$ are consensus and innovation weight sequences given by
\begin{equation}
\label{eq:dist3}
\beta_{t}=\frac{b}{(t+1)^{\delta_{1}}},\alpha_{t}=\frac{a}{t+1},
\end{equation}
where $a, b > 0, 0<\delta_{1}<1/2-1/(2+\epsilon_{1})$ and $\epsilon_{1}$ was defined in Assumption \ref{m:4}.\\
\noindent The update in \eqref{eq:dist1} can be written in a compact manner as follows:
\begin{align}
\label{eq:dist4}
&\wx(t+1)=\mathbf{x}(t)-\beta_{t}\left(\mathbf{L}\otimes\mathbf{I}_{M}\right)\mathbf{x}(t)\nonumber\\&+\alpha_{t}\mathbf{G}(\mathbf{x}(t))\mathbf{R}^{-1}\left(\mathbf{y}(t)-\mathbf{f}\left(\mathbf{x}(t)\right)\right),
\end{align}
where: $\mathbf{x}(t)^{\top}=[\mathbf{x}_{1}(t)^{\top}\cdots \mathbf{x}_{N}(t)^{\top}]$; $\wx(t)^{\top}=[\wx_{1}(t)^{\top}\cdots \wx_{N}(t)^{\top}]$, $\mathbf{f}(\mathbf{x}(t))=\left[\mathbf{f}_{1}(\mathbf{x}_{1}(t))^{\top}\cdots \mathbf{f}_{N}(\mathbf{x}_{N}(t))^{\top}\right]^{\top}$; $\mathbf{R}^{-1}=\mbox{diag}\left[\mathbf{R}_{1}^{-1}, \cdots, \mathbf{R}_{N}^{-1}\right]$; and $\mathbf{G}\left(\mathbf{x}(t)\right)=\mbox{diag}\left[\nabla \mathbf{f}_{1}\left(\mathbf{x}_{1}(t)\right), \cdots, \nabla \mathbf{f}_{N}\left(\mathbf{x}_{N}(t)\right)\right]$.
\noindent We refer to the parameter estimate update in \eqref{eq:dist2} and the projection in \eqref{eq:dist3} as the $\mathcal{CIWNLS}$~($\emph{Consensus}+\emph{innovations}$ Weighted Nonlinear Least Squares) algorithm.
\begin{Remark}
	\label{rm:1}
	The parameter update is recursive and distributed in nature and hence is an \emph{online} algorithm. Moreover, the projection step in \eqref{eq:dist2} ensures that the parameter estimate sequence $\{\mathbf{x}_{n}(t)\}$ is feasible and belongs to the parameter set $\Theta$ at all times $t$.
\end{Remark}
\noindent Methods for analyzing the convergence of distributed stochastic algorithms of the form~\eqref{eq:dist1}-\eqref{eq:dist4} and variants were developed in~\cite{KarMouraRamanan-Est-2008,kar2011convergence,Kar-AdaptiveDistEst-SICON-2012,kar2014asymptotically}. The key is to obtain conditions that ensure the existence of appropriate stochastic Lyapunov functions. To enable this, we propose a condition on the sensing functions (standard in the literature of general recursive procedures) that guarantees the existence of such Lyapunov functions and, hence, the convergence of the distributed estimation procedure.\\

\begin{myassump}{M7}
	\label{m:7}\emph{The following aggregate strict monotonicity condition holds: there exists a constant $c_{1}>0$ such that for each pair $\btheta, \acute{\btheta}$ in $\Theta$ we have that}
\begin{align}
\label{eq:Lyap2}
\sum_{n=1}^{N}\left(\btheta-\acute{\btheta}\right)^{\top} \left(\nabla f_{n}(\btheta)\right)\mathbf{R}_{n}^{-1}\left(f_{n}(\btheta)-f_{n}(\acute{\btheta})\right)\geq c_{1}\left\|\btheta-\acute{\btheta}\right\|^{2}.
\end{align}
\end{myassump}
\noindent In this paper, we assume that the noise covariances are known apriori. However, in scenarios where the noise covariances are not known apriori, in order to verify Assumption \ref{m:7}, only the gradient $\nabla\mathbf{f}_{n}\left(\cdot\right)$ needs to be computed. In case of unknown noise distribution, i.e., unknown noise covariance, the first few observations can be used to estimate the noise covariance so as to get a reasonable estimate of the inverse noise covariance. The estimated noise covariance can then be used to verify the assumption.

\begin{Remark}
	\label{rm:2}
	We comment on the Assumption \ref{m:2}-\ref{m:7}. Assumptions \ref{m:2}-\ref{m:4} are classical with respect to the WNLS convergence.  Assumption M6 specifies some smoothness conditions of the non-linear sensing functions. The smoothness conditions aid in establishing the consistency of the recursive $\mathcal{CIWNLS}$ algorithm. The classical WNLS is usually posed in a non-recursive manner, while the distributed algorithm we propose in this paper is recursive and hence, to ensure convergence we need Lyapunov type conditions, which in turn is specified by Assumption \ref{m:7}. Moreover, Assumptions \ref{m:6}-\ref{m:7} are only sufficient conditions. The key assumptions to establish our main results, Assumptions \ref{m:2}, \ref{m:1}, \ref{m:6}, and \ref{m:7} are required to hold only in the parameter set $\Theta$ and need not hold globally in the entire space $\mathbb{R}^{M}$. This allows our approach to apply to very general nonlinear sensing functions. For example, for functions of the trigonometric type (see Section~\ref{sec:sim} for an illustration), properties such as the strict monotonicity condition in \ref{m:7} hold in the fundamental period, but not globally. As another specific instance, if the $\mathbf{f}_{n}(\cdot)$'s are linear\footnote{To be specific, $\mathbf{f}_{n}\left(\btheta\right)$ is then given by $\mathbf{F}_{n}\btheta$, where $\mathbf{F}_{n}$ is the sensing matrix with dimensions $M_{n}\times M$.}, condition \eqref{WNLS1} in Assumption \ref{m:1}, reduces to $\sum_{n=1}^{N}\mathbf{F}_{n}^{\top}\mathbf{R}_{n}^{-1}\mathbf{F}_{n}$ being full rank (and hence positive definite). The monotonicity condition in Assumption \ref{m:7} in this context coincides with Assumption \ref{m:1}, i.e., it is trivially satisfied by the positive definiteness of the matrix $\sum_{n=1}^{N}\mathbf{F}_{n}^{\top}\mathbf{R}_{n}^{-1}\mathbf{F}_{n}$. Asymptotically efficient distributed parameter estimation schemes for the general linear model have been developed in~\cite{kar2011convergence,Kar-AdaptiveDistEst-SICON-2012}.
\end{Remark}
\section{Main Results}
\label{sec:main_res}
\noindent This section states the main results.
The first concerns the consistency of the estimate sequence in the $\mathcal{CIWNLS}$ algorithm; the proof is in Section \ref{subsec:th_cons}.
\begin{Theorem}
\label{th:cons} Let assumptions~\ref{m:2}-\ref{m:1} and \ref{m:4}-\ref{m:7} hold. Furthermore, assume that the constant $a$ in \eqref{eq:dist3} satisfies
\begin{align}
\label{eq:th:cons0}
ac_{1}\geq 1,
\end{align}
where $c_{1}$ is defined in Assumption \ref{m:7}. Consider the sequence $\{\mathbf{x}_{n}(t)\}$ generated by \eqref{eq:dist2}-\eqref{eq:dist3} at each agent $n$. Then, for each $n$, we have
\begin{equation}
\label{eq:th_cons1}
\mathbb{P}_{\btheta}\left(\lim_{t\rightarrow\infty}(t+1)^{\tau}\|\mathbf{x}_{n}(t)-\btheta\|=0\right)=1
\end{equation}
for all $\tau\in [0, 1/2)$. In particular, the estimate sequence generated by the distributed algorithm~\eqref{eq:dist1}-\eqref{eq:dist4} at any agent $n$ is consistent, i.e., $\mathbf{x}_{n}(t)\rightarrow\btheta$ a.s. as $t\rightarrow\infty$.
%
\end{Theorem}

\noindent At this point, we note that the convergence in Theorem~\ref{th:cons} is order-optimal, in that standard arguments in (centralized) estimation theory show that in general there exists no $\tau\geq 1/2$ such that a centralized WNLS estimator $\{\wbtheta_{t}\}$ satisfies $(t+1)^{\tau}\|\wbtheta_{t}-\btheta\|\rightarrow 0$ a.s. as $t\rightarrow\infty$. \\
\noindent The next result establishes the asymptotic normality of the parameter estimate sequence $\{\mathbf{x}_{n}(t)\}$ and characterizes the asymptotic covariance of the proposed $\mathcal{CIWNLS}$ estimator. This can be benchmarked with the asymptotic covariance of the centralized WNLS estimator. We direct the reader to Section \ref{subsec:th_2} for a proof of Theorem \ref{th:2}.

\begin{Theorem}
	\label{th:2}
	Let the assumptions \ref{m:2}-\ref{m:7} hold. Assume that in addition to assumption \ref{m:2}, the parameter set $\Theta$ is a bounded set. Furthermore, let $a$ defined in \eqref{eq:dist3} satisfy
	\begin{align}
	\label{eq:th2_0}
	a>\max\left\{\frac{1}{c_{1}},\frac{1}{2\inf_{\btheta\in\Theta}\mathbf{\Lambda}_{\btheta,\mbox{\scriptsize{min}}}}\right\},
	\end{align}
	where $c_{1}$ is defined in Assumption \ref{m:7}, $\mathbf{\Lambda}_{\btheta}$ and  $\mathbf{\Lambda}_{\btheta,\mbox{\scriptsize{min}}}$ denote respectively the diagonal matrix of eigenvalues and the minimum eigenvalue of $\mathbf{\Gamma}_{\btheta}$, with $\mathbf{\Gamma}_{\btheta}$ defined in \eqref{eq:m:3_1}.
	Then, for each $n$, the parameter estimate sequence at agent $n$, $\left\{\mathbf{x}_{n}(t)\right\}$, under $\mathbb{P}_{\btheta}$ satisfies the following asymptotic normality condition,
	\begin{align}
	\label{eq:th2_1}
	\sqrt{t+1}\left(\mathbf{x}_{n}(t)-\btheta\right)\overset{\mathcal{D}}{\Longrightarrow}\mathcal{N}\left(0, \mathbf{\Sigma}_{d}\right),
	\end{align}
	
	\noindent where
	\begin{align}
	\label{eq:th2_2}
	\mathbf{\Sigma}_{d}=\frac{\mathbf{aI}}{2N}+\frac{\left(N\mathbf{\Gamma}_{\btheta}-\frac{N\mathbf{I}}{2a}\right)^{-1}}{4},
	\end{align}
	and $\overset{\mathcal{D}}{\Longrightarrow}$ refers to convergence in distribution (weak convergence).
\end{Theorem}

\noindent As the parameter set $\Theta$ is a bounded set in Theorem \ref{th:2}, in addition to being a closed set as in Assumption \ref{m:2}, we have that $\Theta$ is compact and hence  $\inf_{\btheta\in\Theta}\mathbf{\Lambda}_{\btheta,\mbox{\scriptsize{min}}}=\min_{\btheta\in\Theta}\mathbf{\Lambda}_{\btheta,\mbox{\scriptsize{min}}}$, i.e., the infimum is attained. Moreover, from Assumption \ref{m:3}, we have that the matrix $\mathbf{\Gamma}_{\btheta}$ is invertible $\forall \btheta\in\Theta$ and hence $\inf_{\btheta\in\Theta}\mathbf{\Lambda}_{\btheta,\mbox{\scriptsize{min}}} > 0$. Further, under the assumption that $a>\frac{1}{2\inf_{\btheta\in\Theta}\mathbf{\Lambda}_{\btheta,\mbox{\scriptsize{min}}}}$, the difference of the asymptotic covariance of the distributed estimator and that of the centralized estimator, ~i.e., the matrix $\frac{a\mathbf{I}}{2N}+\frac{\left(\mathbf{\Gamma}_{\btheta}-\frac{\mathbf{I}}{2a}\right)^{-1}}{4N}-\left(N\mathbf{\Gamma}_{\btheta}\right)^{-1}$, is positive semidefinite. The above claim can be established by comparing the $i$-th eigenvalue of the asymptotic covariance of the distributed estimator with that of the centralized estimator, as both the covariance matrices are simultaneously diagonalizable, i.e., have the same set of eigenvectors. To be specific,
\begin{align}
\label{eq:asymp_var_2}
&a^{2}\mathbf{\Lambda}_{\btheta,ii}^{2}-2a\mathbf{\Lambda}_{\btheta,ii}+1 \geq 0\nonumber\\
&\Rightarrow\frac{1}{\mathbf{\Lambda}_{\btheta,ii}}\le \frac{a^{2}\mathbf{\Lambda}_{\btheta,ii}}{2a\mathbf{\Lambda}_{\btheta,ii}-1}\nonumber\\
&\Rightarrow\frac{1}{N\mathbf{\Lambda}_{\btheta,ii}}\le \frac{a^{2}\mathbf{\Lambda}_{\btheta,ii}}{2aN\mathbf{\Lambda}_{\btheta,ii}-N},
\end{align}
which holds for all $i=1,2,\cdots, N$.\\
\noindent We now benchmark the asymptotic covariance of the proposed estimator $\mathcal{CIWNLS}$ with that of the optimal centralized estimator. From Assumption \ref{m:6}, we have for all $\btheta\in\Theta$
\begin{align}
\label{eq:asymp_var_3}
\left\|\mathbf{\Gamma}_{\btheta}\right\|\le\max_{n=1,\cdots, N}k_{n}^{2}\left\|\mathbf{R}_{n}^{-1}\right\|=k_{\mbox{\scriptsize{max}}}^{*},
\end{align}
where $k_{n}$ is defined in Assumption \ref{m:6}. Moreover, from the hypotheses of Theorem \ref{th:2} we have that $\mathbf{\Lambda}_{\btheta,\mbox{\scriptsize{min}}}>\frac{1}{2a}$, for all $\btheta\in\Theta$. Thus, we have the following characterization of the eigenvalues for the matrix $\mathbf{\Gamma}_{\btheta}$ for all $\btheta\in\Theta$,
\begin{align}
\label{eq:asymp_var_4}
\frac{1}{2a} < \mathbf{\Lambda}_{\btheta,ii} \le k_{\mbox{\scriptsize{max}}}^{*},
\end{align}
for all $i$. The difference of the $i$-th eigenvalue of the asymptotic covariance of the distributed estimator and the centralized estimator $\mathbf{\Lambda}_{d,i}$, is given by $\mathbf{\Lambda}_{d,i}=\frac{\left(a\mathbf{\Lambda}_{\btheta,ii}-1\right)^{2}}{N\mathbf{\Lambda}_{\btheta,ii}\left(2a\mathbf{\Lambda}_{\btheta,ii}-1\right)}$.
Now, we consider two cases. Specifically, if the condition
\begin{align}
\label{eq:asymp_var_41}
k_{\mbox{\scriptsize{max}}}^{*} > \max\left\{\frac{1}{c_{1}},\frac{1}{2\inf_{\btheta\in\Theta}\mathbf{\Lambda}_{\btheta,\mbox{\scriptsize{min}}}}\right\},
\end{align}
is satisfied, then $a$ can be chosen to be $a<1/k_{\mbox{\scriptsize{min}}}^{*}$, and then we have,
\begin{align}
\label{eq:asymp_var_5}
\frac{1}{2a} < \mathbf{\Lambda}_{\btheta,ii} \le \frac{1}{a}.
\end{align}
It is to be noted that the function $h(x)=\frac{\left(ax-1\right)^{2}}{Nx\left(2ax-1\right)}$ is non-increasing in the interval $\left(\frac{1}{2a}, \frac{1}{a}\right)$. Hence, we have that
\begin{align}
\label{eq:asymp_var_6}
\left\|\mathbf{\Sigma}_{d}-\mathbf{\Sigma}_{c}\right\|= \frac{\left(a\mathbf{\Lambda}_{\btheta,\mbox{\scriptsize{min}}}-1\right)^{2}}{N\mathbf{\Lambda}_{\btheta,\mbox{\scriptsize{min}}}\left(2a\mathbf{\Lambda}_{\btheta,\mbox{\scriptsize{min}}}-1\right)},
\end{align}
where,
\begin{align}
\label{eq:asymp_var_7}
\max\left\{\frac{1}{c_{1}},\frac{1}{2\inf_{\btheta\in\Theta}\mathbf{\Lambda}_{\btheta,\mbox{\scriptsize{min}}}}\right\}<a<\frac{1}{k_{\mbox{\scriptsize{max}}}^{*}}.
\end{align}
In the case, when the condition in \eqref{eq:asymp_var_41} is violated, we have that,
{\small\begin{align}
\label{eq:asymp_var_8}
&\left\|\mathbf{\Sigma}_{d}-\mathbf{\Sigma}_{c}\right\|\nonumber\\ &=\max\left\{\frac{\left(a\mathbf{\Lambda}_{\btheta,\mbox{\scriptsize{min}}}-1\right)^{2}}{N\mathbf{\Lambda}_{\btheta,\mbox{\scriptsize{min}}}\left(2a\mathbf{\Lambda}_{\btheta,\mbox{\scriptsize{min}}}-1\right)}, \frac{\left(a\mathbf{\Lambda}_{\btheta,\mbox{\scriptsize{max}}}-1\right)^{2}}{N\mathbf{\Lambda}_{\btheta,\mbox{\scriptsize{max}}}\left(2a\mathbf{\Lambda}_{\btheta,\mbox{\scriptsize{max}}}-1\right)}\right\}\nonumber\\
& \le \max\left\{\frac{\left(a\mathbf{\Lambda}_{\btheta,\mbox{\scriptsize{min}}}-1\right)^{2}}{N\mathbf{\Lambda}_{\btheta,\mbox{\scriptsize{min}}}\left(2a\mathbf{\Lambda}_{\btheta,\mbox{\scriptsize{min}}}-1\right)}, \frac{\left(ak_{\mbox{\scriptsize{max}}}^{*}-1\right)^{2}}{Nk_{\mbox{\scriptsize{max}}}^{*}\left(2ak_{\mbox{\scriptsize{max}}}^{*}-1\right)}\right\},
\end{align}}
where $\mathbf{\Lambda}_{\btheta,\mbox{\scriptsize{max}}}$ denotes the largest eigenvalue of $\mathbf{\Gamma}_{\btheta}$. Note that the proposition in \eqref{eq:asymp_var_8} is equivalent to \eqref{eq:asymp_var_6}, when the condition in \eqref{eq:asymp_var_41} is satisfied. Hence, for all feasible choices of $a$, which are in turn given by \eqref{eq:th2_0}, the characterization in \eqref{eq:asymp_var_8} holds.\\
\noindent The above mentioned findings can be precisely stated in the form of the following corollary:
\begin{corollary}
	\label{cor:1}
	\noindent Let the hypothesis of Theorem \ref{th:2} hold. Then, we have,
	{\small\begin{align}
	\label{eq:cor1}
	&\left\|\mathbf{\Sigma}_{d}-\mathbf{\Sigma}_{c}\right\|\nonumber\\ 
	& \le \max\left\{\frac{\left(a\mathbf{\Lambda}_{\btheta,\mbox{\scriptsize{min}}}-1\right)^{2}}{N\mathbf{\Lambda}_{\btheta,\mbox{\scriptsize{min}}}\left(2a\mathbf{\Lambda}_{\btheta,\mbox{\scriptsize{min}}}-1\right)}, \frac{\left(ak_{\mbox{\scriptsize{max}}}^{*}-1\right)^{2}}{Nk_{\mbox{\scriptsize{max}}}^{*}\left(2ak_{\mbox{\scriptsize{max}}}^{*}-1\right)}\right\},
	\end{align}}
\noindent where $\mathbf{\Sigma}_{d}$ and $\mathbf{\Sigma}_{c}$ are defined in \eqref{eq:th2_2} and \eqref{eq:asymp_var_11}, respectively.
\end{corollary}
\noindent Furthermore, as noted above that the difference of the asymptotic covariance of the distributed estimator and that of the centralized estimator is positive semi-definite.~(This is intuitively expected, as a distributed procedure may not outperform a well-designed centralized procedure.) The inefficiency of the distributed estimator with respect to the centralized WNLS estimator, as far as the asymptotic covariance is concerned, is due to the use of suboptimal innovation gains (see, for example \eqref{eq:dist1}) used in the parameter estimate update. An optimal innovation gain sequence would require the knowledge of the global model information, i.e., the sensing functions and the noise covariances across all the agents.~(See Remark 4.3 below for a detailed discussion.) Though the distributed estimation scheme is suboptimal with respect to the centralized estimator as far as the asymptotic covariance is concerned, its performance is significantly better than the non-collaborative case, i.e., in which agents perform estimation in a fully decentralized or isolated manner. In particular, since the agent sensing models are likely to be locally unobservable for $\btheta$, the asymptotic covariances in the non-collaborative scenario may diverge due to non-degeneracy.
\begin{Remark}
	\label{rm:3}
	In this context, we briefly review the methodology adopted in \cite{kar2014asymptotically} to achieve asymptotically efficient distributed estimators for general standard statistical exponential families. The asymptotic efficiency of the estimator proposed in \mbox{\cite{kar2014asymptotically}} is a result of a certainty-equivalence type distributed optimal gain sequence generated through an auxiliary consistent parameter estimate sequence (see, for example Section III.A. in \mbox{\cite{kar2014asymptotically}}). The generation of the auxiliary estimate sequence comes at the cost of more communication and computation complexity as two other parallel recursions run in addition to the parameter estimate recursion. To be specific, the adaptive gain refinement, which is the key for achieving asymptotic efficiency, involves communication of the gain matrices that belong to the space $\mathbb{R}^{M\times M}$. Moreover, in a setting like the one considered in this paper, where the parameter belongs to a closed convex subset $\Theta\in\mathbb{R}^{M}$, the parameter estimate at all times provided by the algorithm in \mbox{\cite{kar2014asymptotically}} might not be feasible. In contrast, the communication and computation complexity in $\mathcal{CIWNLS}$ is significantly lower than that of the algorithm proposed in \mbox{\cite{kar2014asymptotically}}. The price paid by $\mathcal{CIWNLS}$ is the lower asymptotic performance as measured in terms of the asymptotic covariance as we discuss next.
\end{Remark}
\noindent We compare the computational and communication overhead of $\mathcal{CIWNLS}$ and of the algorithm in \cite{kar2014asymptotically}. For simplicity, we consider a $d$-regular communication graph with every agent connected to $d$ other agents. We compute the computation and communication overhead agent wise. In one sampling epoch, an agent in $\mathcal{CIWNLS}$ communicates $M$-dimensional parameter estimates to its neighbors, i.e., the communication overhead is $Md$. In the algorithm proposed in \cite{kar2014asymptotically}~(see, (8)-(11) in Section $III.A$), an agent not only communicates its auxiliary and optimal parameter estimates but also its gain matrix, to its neighbors, with the communication overhead $2Md+M^{2}d$. With respect to the computational overhead, in every sampling epoch the number of computations in $\mathcal{CIWNLS}$ at agent $n$ is given by $O\left(M_{n}M+M(d+1)\right)$. The maximum computational overhead across all agents is thus given by $\max_{n=1,\cdots,N}O\left(M_{n}M+M(d+1)\right)$. In comparison, the number of computations at any agent in the algorithm proposed in \cite{kar2014asymptotically} is given by $O\left(M^{3}+2M^{2}(d+1)+2M(d+1)\right)$. Thus, the communication and computational complexity of the proposed algorithm is much lower than that of the algorithm proposed in \cite{kar2014asymptotically}, at the cost of suboptimal asymptotic estimation error covariances.
\section{Simulations}
\label{sec:sim}
\begin{figure}
	\centering
	\captionsetup{justification=centering}
	\includegraphics[width=90mm]{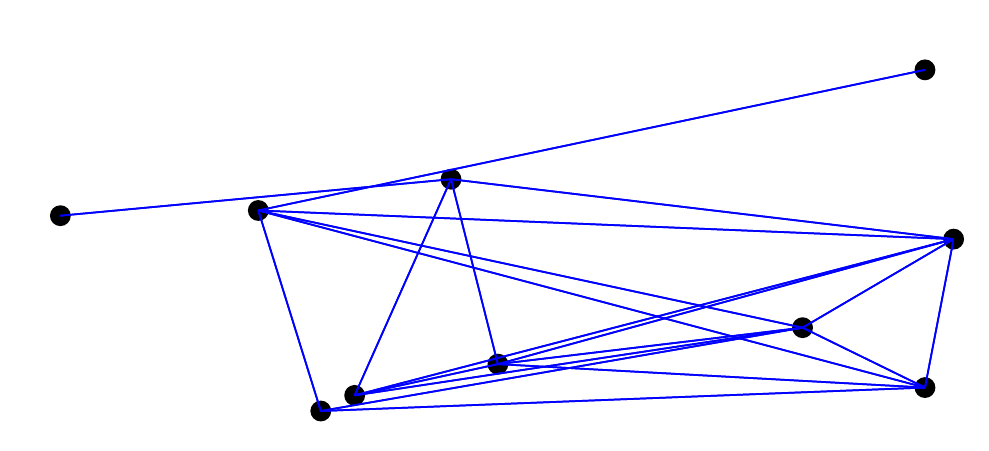}
	\caption{Network Deployment of $10$ agents}\label{fig:1}
\end{figure}

\begin{figure}
	\centering
	\captionsetup{justification=centering}
	\includegraphics[width=90mm]{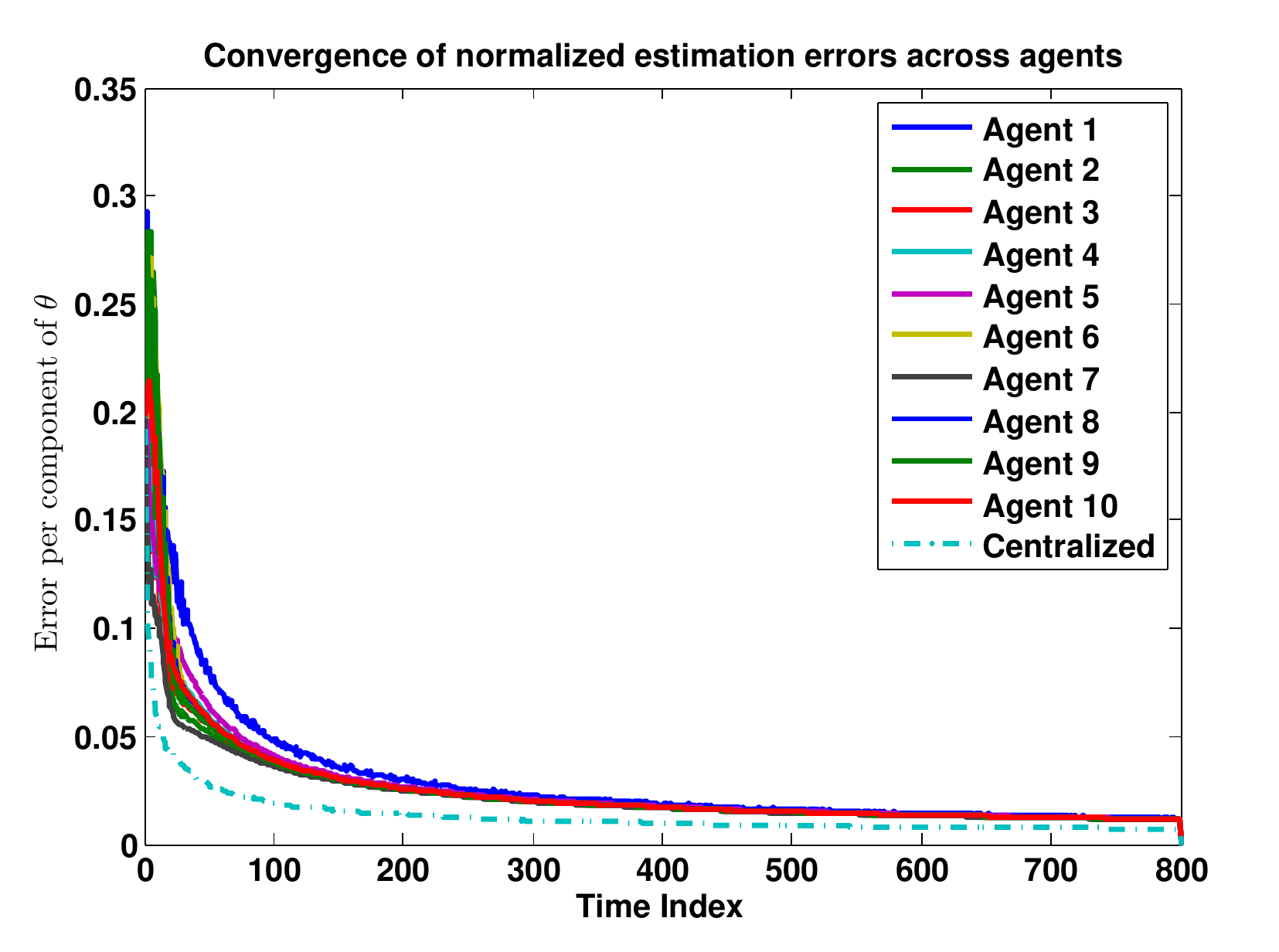}
	\caption{Convergence of normalized estimation error at each agent}\label{fig:2}
\end{figure}

\noindent We generate a random geometric network of $10$ agents, shown in Figure \ref{fig:1}.
The $x$ coordinates and the $y$ coordinates of the agents are sampled from a uniform distribution on the interval $[0,1]$. We link two vertices by an edge if the distance between them is less than or equal to $g=0.4$. We go on re-iterating this procedure until we get a connected graph. We choose the parameter set $\Theta$ to be $\Theta=\left[-\frac{\pi}{4}, \frac{\pi}{4}\right]^{5}\in\mathbb{R}^{5}$. This choice of $\Theta$ conforms with Assumption \ref{m:2}. The sensing functions are chosen to be certain trigonometric functions as described below. The underlying parameter is $5$ dimensional, $\btheta=\left[\theta_{1},~\theta_{2},~\theta_{3},~\theta_{4},~\theta_{5}\right]$. The sensing functions across different agents are given by, $\mathbf{f}_{1}(\btheta)=\sin(\theta_{1}+\theta_{2}), \mathbf{f}_{2}(\btheta)=\sin(\theta_{3}+\theta_{2}), \mathbf{f}_{3}(\btheta)=\sin(\theta_{3}+\theta_{4}), \mathbf{f}_{4}(\btheta)=\sin(\theta_{4}+\theta_{5}), \mathbf{f}_{5}(\btheta)=\sin(\theta_{1}+\theta_{5}), \mathbf{f}_{6}(\btheta)=\sin(\theta_{1}+\theta_{3}), \mathbf{f}_{7}(\btheta)=\sin(\theta_{4}+\theta_{2}), \mathbf{f}_{8}(\btheta)=\sin(\theta_{3}+\theta_{5}), \mathbf{f}_{9}(\btheta)=\sin(\theta_{1}+\theta_{4})$ and  $\mathbf{f}_{10}(\btheta)=\sin(\theta_{1}+\theta_{5})$. Clearly, the local sensing models are unobservable, but collectively they are globally observable since, in the parameter set $\Theta$ under consideration, $\sin(\cdot)$ is one-to-one and the set of linear combinations of the $\btheta$ components corresponding to the arguments of the $\sin(\cdot)$'s constitute a full-rank system for $\btheta$. Hence, the sensing model conforms to Assumption \ref{m:1}. The agents make noisy scalar observations where the observation noise process is Gaussian and the noise covariance is given by $\mathbf{R}=2\mathbf{I}_{10}$. The true (but unknown) value of the parameter is taken to be $\btheta=\left[\pi/6,~-\pi/7,~\pi/12,~-\pi/5,~\pi/16\right]$. It is readily verified that this sensing model and the parameter set $\Theta=\left[-\frac{\pi}{4}, \frac{\pi}{4}\right]^{5}$ satisfy Assumptions \ref{m:3}-\ref{m:7}. The projection operator $\mathcal{P}_{\Theta}$ onto the set $\Theta$ defined in \eqref{eq:dist2} is given by,
\begin{align}
\label{eq:sim_1}
\left[\mathbf{x}_{n}(t)\right]_{i}=
\begin{cases}
\frac{\pi}{4} & [\wx_{n}(t)]_{i} \geq \frac{\pi}{4}\\
\left[\wx_{n}(t)\right]_{i} & \frac{-\pi}{4} < [\wx_{n}(t)]_{i} < \frac{\pi}{4}\\
\frac{-\pi}{4} & [\wx_{n}(t)]_{i} < \frac{-\pi}{4},
\end{cases}
\end{align}
for all $i=1,\cdots, M$. \\
The sensing model is motivated by distributed static phase estimation in smartgrids. For a more complete treatment of the classical problem of static phase estimation in power grids, we direct the reader to \cite{Ilic-Zaborsky}. Coming back to the current context, to be specific, the physical grid can be modeled as a network with the loads and generators being the nodes~(vertices), while the transmission lines being the edges, and the sensing model reflects the power flow equations. The goal of distributed static phase estimation is to estimate the vector of phases from line flow data. The interested reader is directed to Section $IV.D$ of \cite{KarMouraRamanan-Est-2008} for a detailed treatment of distributed static phase estimation.\\
We carry out $250$ Monte-Carlo simulations for analyzing the convergence of the parameter estimates and their asymptotic covariances. The estimates are initialized to be $0$, i.e., $\mathbf{x}_{n}(0)=\mathbf{0}$ for $n=1,\cdots, 5$. The normalized error for the $n$-th agent at time $t$ is given by the quantity $\left\|\mathbf{x}_{n}(t)-\btheta\right\|/5$. Figure \ref{fig:2} shows the normalized error at every agent against the time index $t$. We compare it with the normalized error of the centralized estimator in Figure \ref{fig:2}.
\begin{figure}
	\centering
	\captionsetup{justification=centering}
	\includegraphics[width=90mm]{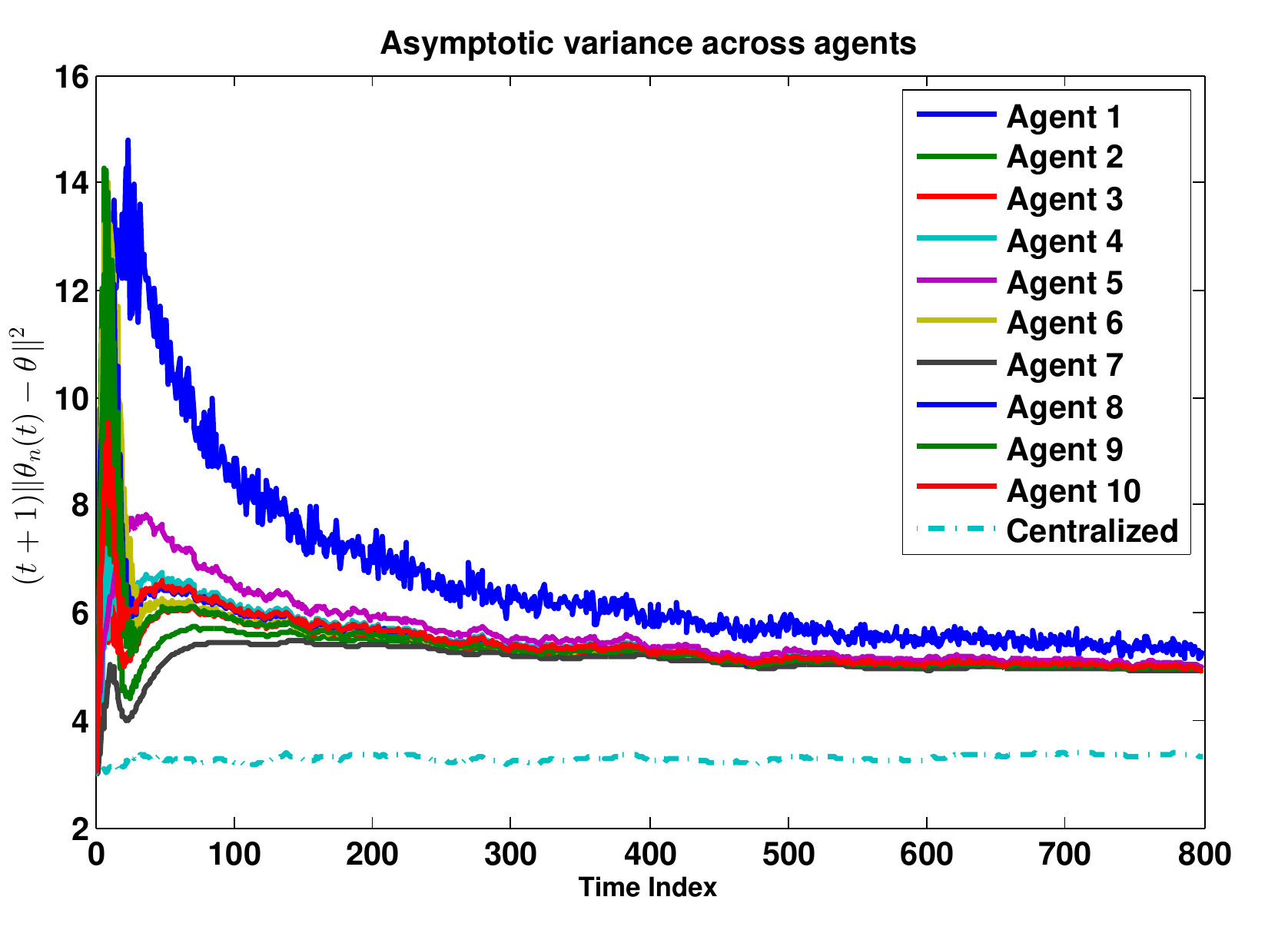}
	\caption{Asymptotic variance at each agent}\label{fig:3}
\end{figure}
We note that the errors converge to zero as established in Theorem \ref{th:cons}. The decrease in error is rapid in the beginning and slows down with increasing $t$; this is a consequence of the decreasing weight sequences $\{\alpha_{t}\}$ and $\{\beta_{t}\}$.
Finally, in Fig. \ref{fig:3} we compare the asymptotic variances of our scheme and that of the centralized WNLS estimator. For the distributed estimator $\mathcal{CIWNLS}$ and for each $n$, Fig. \ref{fig:3} plots the quantities {\small$(t+1)\left\|\mathbf{x}_{n}(t)-\btheta\right\|^{2}$} averaged over the Monte-Carlo trials.
By Theorem 4.2, this quantity is expected to converge to the trace of
the asymptotic covariance $\mathbf{\Sigma}_{d}$ of the $\mathcal{CIWNLS}$, i.e., $\mbox{tr}(\mathbf{\Sigma}_{d})$, the same for all $n$. We also simulate the centralized WNLS and plot the scaled error {\small$(t+1)\left\|\wbtheta(t)-\btheta\right\|^{2}$} averaged over the Monte-Carlo trials. Similarly, from Proposition \ref{prop:WNLS} we have that {\small$\mathbb{E}_{\btheta}\left[(t+1)\left\|\wbtheta(t)-\btheta\right\|^{2}\right]\to \mbox{tr}\left(\mathbf{\Sigma}_{c}\right)$}. In this simulation setup, $\mbox{tr}\left(\mathbf{\Sigma}_{c}\right)$ and $\mbox{tr}\left(\mathbf{\Sigma}_{d}\right)$ are evaluated to be $3.6361$ and $5.4517$, respectively, a loss of about 1.76 dB. From the simulation experiment conducted above, the experimental values of $\mbox{tr}\left(\mathbf{\Sigma}_{c}\right)$ and $\mbox{tr}\left(\mathbf{\Sigma}_{d}\right)$ are found to be $3.9554$ and $5.6790$ respectively. These experimental findings verify the conclusions of Theorem 4.2.
\section{Proof of Main Results}
\label{sec:proof_main_res}
In this section, we provide the proofs of Theorems \ref{th:cons} and \ref{th:2}.
\subsection{Proof of Theorem \ref{th:cons}}
\label{subsec:th_cons}
\begin{IEEEproof}
	The proof of Theorem \ref{th:cons} is accomplished in three steps. First, we establish the boundedness of the estimate sequence followed by proving the strong consistency of the estimate sequence $\{\mathbf{x}_{n}(t)\}$ and then in the sequel we establish the rate of convergence of the estimate sequence to the true underlying parameter. We follow the basic idea developed in \cite{kar2014asymptotically}.
	\begin{Lemma}
		\label{le:l0}
		Let the hypothesis of Theorem \ref{th:cons} hold. Then, for each $n$, the process $\{\mathbf{x}_{n}(t)\}$ satisfies
		\label{l2}
		{\small\begin{align}
		\label{eq:7.1}
		\mathbb{P}_{\theta}\left(\sup_{t\ge 0} \left\|\mathbf{x}(t)\right\| < \infty\right) =1.
		\end{align}}
	\end{Lemma}
	\begin{IEEEproof}
		\noindent The proof is built around a similar framework as the proof of Lemma IV.1 in \cite{kar2014asymptotically} with appropriate modifications to take into account the state-dependent nature of the innovation gain and the projection operator used in \eqref{eq:dist3}. Define the process $\{\mathbf{z}(t)\}$ as follows $\mathbf{z}(t)=\mathbf{x}(t)-\mathbf{1}_{N} \otimes \btheta$ where $\btheta$ denotes the true~(but unknown) parameter.
		Note the following recursive relationship:
		{\small\begin{align}
		\label{eq:l0_pr_1}
		& \wx(t+1)-\mathbf{1}_{N} \otimes \btheta=\mathbf{z}(t)-\beta_{t}(\mathbf{L}\otimes \mathbf{I}_{M})z(t)\nonumber\\&+\alpha_{t}\mathbf{G}\left(\mathbf{x}(t)\right)\mathbf{R}^{-1}\left(\mathbf{y}(t)-\mathbf{f}(\mathbf{x}(t))\right),
		\end{align}}
		which further implies that,
		{\small\begin{align}
		\label{eq:l0_pr1_1}
		&\wx(t+1)-\mathbf{1}_{N} \otimes \btheta=\mathbf{z}(t)-\beta_{t}(\mathbf{L} \otimes \mathbf{I}_{M})z(t)\nonumber\\&+\alpha_{t}\mathbf{G}\left(\mathbf{x}(t)\right)\left(\mathbf{y}(t)-\mathbf{f}\left(\mathbf{1}_{N}\otimes\btheta\right)\right)\nonumber\\&-\alpha_{t}\mathbf{G}\left(\mathbf{x}(t)\right)\mathbf{R}^{-1}\left(\mathbf{f}\left(\mathbf{x}(t)\right)-\mathbf{f}\left(\mathbf{1}_{N}\otimes\btheta\right)\right).
		\end{align}}
		\noindent In the above, we have used a basic property of the Laplacian $\mathbf{L}$,
		{\small\begin{align}
		\label{eq:l0_pr_2}
		\left(\mathbf{L} \otimes \mathbf{I}_{M}\right)\left(\mathbf{1}_{N} \otimes \btheta\right)=\mathbf{0},
		\end{align}}
		\noindent Since the projection is onto a convex set it is non-expansive. It follows that the inequality
		{\small\begin{align}
		\label{eq:l0_pr_31}
		\left\|\mathbf{x}_{n}(t+1)-\btheta\right\|\leq \left\|\wx_{n}(t+1)-\btheta\right\|
		\end{align}}
		holds for all $n$ and $t$.
		\noindent Taking norms on both sides of \eqref{eq:l0_pr_1} and using \eqref{eq:l0_pr_31}, we have,
		{\small\begin{align}
		\label{eq:l0_pr_3}
		&\left\|\mathbf{z}(t+1)\right\|^{2}\le\left\|\mathbf{z}(t)\right\|^{2}-2\beta_{t}\mathbf{z}^{\top}(t)\left(\mathbf{L} \otimes \mathbf{I}_{M}\right)\mathbf{z}(t) \nonumber\\&-2\alpha_{t}\mathbf{z}^{\top}(t)\mathbf{G}\left(\mathbf{x}(t)\right)\mathbf{R}^{-1}\left(\mathbf{f}\left(\mathbf{x}(t)\right)-\mathbf{f}\left(\mathbf{1}_{N}\otimes\btheta\right)\right)\nonumber\\
		&+\beta_{t}^{2}\mathbf{z}^{\top}(t)\left(\mathbf{L}\otimes \mathbf{I}_{M}\right)^{2}\mathbf{z}(t)\nonumber\\&+2\alpha_{t}\beta_{t}\mathbf{z}^{\top}(t)\left(\mathbf{L}\otimes \mathbf{I}_{M}\right)\mathbf{G}\left(\mathbf{x}(t)\right)\mathbf{R}^{-1}\left(\mathbf{f}\left(\mathbf{x}(t)\right)-\mathbf{f}\left(\mathbf{1}_{N}\otimes\theta\right)\right)\nonumber\\
		&+\alpha_{t}^{2}\left\|\mathbf{G}\left(\mathbf{x}(t)\right)\mathbf{R}^{-1}\left(\mathbf{y}(t)-\mathbf{f}\left(\mathbf{1}_{N}\otimes\btheta\right)\right)\right\|^{2}\nonumber\\&+\alpha_{t}^{2}\left\|\mathbf{G}\left(\mathbf{x}(t)\right)\mathbf{R}^{-1}\left(\mathbf{f}\left(\mathbf{x}(t)\right)-\mathbf{f}\left(\mathbf{1}_{N}\otimes\btheta\right)\right)\right\|^{2}\nonumber\\
		&+2\alpha_{t}\mathbf{z}^{\top}(t)\mathbf{G}\left(\mathbf{x}(t)\right)\mathbf{R}^{-1}\left(\mathbf{y}(t)-\mathbf{f}(\mathbf{1}_{N}\otimes\btheta)\right)\nonumber\\&+2\alpha_{t}^{2}\left(\mathbf{y}(t)-\mathbf{f}(\mathbf{1}_{N}\otimes\btheta)\right)^{\top}\mathbf{R}^{-1}\mathbf{G}^{\top}\left(\mathbf{x}(t)\right)\nonumber\\&\times\mathbf{G}\left(\mathbf{x}(t)\right)\mathbf{R}^{-1}\left(\mathbf{f}\left(\mathbf{1}_{N}\otimes\btheta\right)-\mathbf{f}\left(\mathbf{x}(t)\right)\right).
		\end{align}}
		\noindent Consider the orthogonal decomposition
		{\small\begin{align}
		\label{eq:l0_pr_4}
		\mathbf{z}=\mathbf{z}_{c}+\mathbf{z}_{c \perp},
		\end{align}}
		\noindent where $\mathbf{z}_{c}$ denotes the projection of $\mathbf{z}$ to the consensus subspace {\small$\mathcal{C}=\{\mathbf{z} \in \mathbb{R}^{MN} |\mathbf{z}=1_{N}\otimes a, \mbox{for~~some~~a} \in \mathbb{R}^{M} \}$}.
		\noindent From \eqref{eq:obsmod}, we have that,
		{\small\begin{align}
		\label{eq:l0_pr_6}
		\mathbb{E}_{\btheta}\left[\mathbf{y}(t)-\mathbf{f}\left(\mathbf{1}_{N}\otimes\btheta\right)\right]=\mathbf{0}.
		\end{align}}
		\noindent Consider the process
		{\small\begin{align}
		\label{eq:l0_pr_7}
		V_{2}(t)=\left\|\mathbf{z}(t)\right\|^{2}.
		\end{align}}
		\noindent Using conditional independence properties, we have,
		{\small\begin{align}
		\label{eq:l0_pr_8}
		&\mathbb{E}_{\btheta}[V_{2}(t+1)|\mathcal{F}_{t}]\le V_{2}(t) + \beta_{t}^{2}\mathbf{z}^{\top}(t)\left(\mathbf{\overline{L}} \otimes \mathbf{I}_{M}\right)^{2}\mathbf{z}(t)\nonumber\\
		&+\alpha_{t}^{2}\mathbb{E}_{\btheta}\left[\left\|\mathbf{G}\left(\mathbf{x}(t)\right)\mathbf{R}^{-1}\left(\mathbf{y}(t)-\mathbf{f}\left(\mathbf{1}_{N}\otimes\btheta\right)\right)\right\|^{2}\right]\nonumber\\&-2\beta_{t}\mathbf{z}^{\top}(t)\left(\mathbf{\overline{L}}\otimes \mathbf{I}_{M}\right)\mathbf{z}(t)\nonumber\\
		&-2\alpha_{t}\mathbf{z}^{\top}(t)\mathbf{G}\left(\mathbf{x}(t)\right)\mathbf{R}^{-1}\left(\mathbf{f}\left(\mathbf{x}(t)\right)-\mathbf{f}\left(\mathbf{1}_{N}\otimes\btheta\right)\right)\nonumber\\&+2\alpha_{t}\beta_{t}\mathbf{z}^{\top}(t)\left(\mathbf{\overline{L}} \otimes \mathbf{I}_{M}\right)\mathbf{G}\left(\mathbf{x}(t)\right)\mathbf{R}^{-1}\left(\mathbf{f}\left(\mathbf{x}(t)\right)-\mathbf{f}\left(\mathbf{1}_{N}\otimes\btheta\right)\right) \nonumber\\
		&+\alpha_{t}^{2}\left\|\left(\mathbf{f}\left(\mathbf{x}(t)\right)-\mathbf{f}\left(\mathbf{1}_{N}\otimes\btheta\right)\right)^{\top}\mathbf{G}^{\top}\left(\mathbf{x}(t)\right)\mathbf{R}^{-1}\right\|^{2}.
		\end{align}}
		\noindent We use the following inequalities $\forall t \ge t_{1}$,
		{\small\begin{align}
		\label{eq:l0_pr_9}
		&\mathbf{z}^{\top}(t)\left(\mathbf{L} \otimes \mathbf{I}_{M}\right)^{2}\mathbf{z}(t) \overset{(q1)}{\le} \lambda_{N}^{2}(\mathbf{L})||\mathbf{z}_{C\perp}(t)||^{2};\nonumber\\
		&\mathbf{z}^{\top}(t)\mathbf{G}\left(\mathbf{x}(t)\right)\mathbf{R}^{-1}\left(\mathbf{f}\left(\mathbf{x}(t)\right)-\mathbf{f}\left(\mathbf{1}_{N}\otimes\btheta\right)\right)\ge c_{1}||\mathbf{z}(t)||^{2} \overset{(q2)}{\ge} 0;\nonumber\\
		&\mathbf{z}^{\top}(t)\left(\mathbf{L} \otimes \mathbf{I}_{M}\right)\mathbf{z}(t) \overset{(q3)}{\ge} \lambda_{2}(\mathbf{\overline{L}})\left\|\mathbf{z}_{C\perp}(t)\right\|^{2};\nonumber\\
		&\mathbf{z}^{\top}(t)\left(\mathbf{L} \otimes \mathbf{I}_{M}\right)\mathbf{G}\left(\mathbf{x}(t)\right)\mathbf{R}^{-1}\left(\mathbf{f}\left(\mathbf{x}(t)\right)-\mathbf{f}\left(\mathbf{1}_{N}\otimes\btheta\right)\right) \nonumber\\&\overset{(q4)}{\le} c_{2}\left\|\mathbf{z}(t)\right\|^{2},
		\end{align}}
		\noindent for $c_{1}$ as defined in Assumption \ref{m:7} and a positive constant $c_{2}$. Inequalities $(q1)$ and $(q4)$ follow from the properties of the Laplacian. Inequality $(q2)$ follows from Assumption \ref{m:7} and $(q4)$ follows from Assumption \ref{m:6} since we have that {\small$\left\|\nabla \mathbf{f}_{n}\left(\mathbf{x}_{n}(t)\right)\right\|$} is uniformly bounded from above by $k_{n}$ for all $n$ and hence, we have that {\small$\left\|\mathbf{G}\left(\mathbf{x}(t)\right)\right\|\le \max_{n=1,\cdots,N}k_{n}$}.
		\noindent We also have
		{\small\begin{align}
		\label{eq:l0_pr_13}
		\mathbb{E}_{\btheta}\left[\left\|\mathbf{G}\left(\mathbf{x}(t)\right)\mathbf{R}^{-1}\left(\mathbf{y}(t)-\mathbf{f}\left(\mathbf{1}_{N}\otimes\btheta\right)\right)\right\|^{2}\right] \le c_{4},
		\end{align}}
		\noindent for some constant $c_{4} > 0$. In \eqref{eq:l0_pr_13}, we use the fact that the noise process under consideration has finite covariance. We also use the fact that {\small$\left\|\mathbf{G}\left(\mathbf{x}(t)\right)\right\|\le \max_{n=1,\cdots,N}k_{n}$}, which in turn follows from Assumption \ref{m:5}.
		\noindent We further have that,
		{\small\begin{align}
		\label{eq:l0_pr_14}
		\left\|\mathbf{G}\left(\mathbf{x}(t)\right)\mathbf{R}^{-1}\left(\mathbf{f}\left(\mathbf{x}(t)\right)-\mathbf{f}\left(\mathbf{1}_{N}\otimes\btheta\right)\right)\right\|^{2}\le c_{3}\left\|\mathbf{z}(t)\right\|^{2},
		\end{align}}
		\noindent where $c_{3} > 0$ is a constant.
		\noindent It is to be noted that \eqref{eq:l0_pr_14} follows from the Lipschitz continuity in Assumption \ref{m:5} and the result that {\small$\left\|\mathbf{G}\left(\mathbf{x}(t)\right)\right\|\le \max_{n=1,\cdots,N}k_{n}$}.
		\noindent Using \eqref{eq:l0_pr_8}-\eqref{eq:l0_pr_14}, we have,
		{\small\begin{align}
		\label{eq:l0_pr_15}
		&\mathbb{E}_{\btheta}[V_{2}(t+1)|\mathcal{F}_{t}]\le \left(1+c_{5}\left(\alpha_{t}\beta_{t}+\alpha_{t}^{2}\right)\right)V_{2}(t)\nonumber\\&-c_{6}(\beta_{t}-\beta^{2}_{t})||\mathbf{x}_{C\perp}(t)||^{2}+c_{4}\alpha_{t}^{2},
		\end{align}}
		\noindent for some positive constants $c_{5}$ and $c_{6}$.
		\noindent As $\beta_{t}^{2}$ goes to zero faster than $\beta_{t}$, $\exists t_{2}$ such that $\forall t\ge t_{2}$, $\beta_{t} \ge \beta^{2}_{t}$.
		\noindent Hence $\exists t_{2}$ and $\exists \tau_{1}, \tau_{2} > 1$  such that for all $t \ge t_{2}$
		{\small\begin{align}
		\label{eq:l0_pr_16}
		c_{5}\left(\alpha_{t}\beta_{t}+\alpha_{t}^{2}\right) \le \frac{c_{7}}{(t+1)^{\tau_{1}}}=\gamma_{t}~~\mbox{and}~~
		c_{4}\alpha_{t}^{2}\le\frac{c_{8}}{(t+1)^{\tau_{2}}}=\hat{\gamma}_{t},
		\end{align}}
		\noindent where $c_{7}, c_{8} > 0$ are constants.
		\noindent By the above construction we obtain $\forall t \geq t_{2}$,
		{\small\begin{align}
		\label{eq:l2_pr_17}
		\mathbb{E}_{\theta^{*}}[V_{2}(t+1) | \mathcal{F}_{t}] \le (1+\gamma_{t})V_{2}(t)+\hat{\gamma}_{t},
		\end{align}}
		\noindent where the positive weight sequences $\{\gamma_{t}\}$ and $\{\hat{\gamma}_{t}\}$ are summable i.e.,
		{\small\begin{align}
		\label{eq:l2_pr_18}
		\sum_{t\ge0}\gamma_{t} < \infty, 	\sum_{t\ge0}\hat{\gamma}_{t} < \infty.
		\end{align}}
		\noindent By \eqref{eq:l2_pr_18}, the product $\prod_{s=t}^{\infty}(1+\gamma_{s})$ exists for all $t$. Now let $\{W(t)\}$ be such that
		{\small\begin{align}
		\label{eq:l2_pr_19}
		W(t)=\left(\prod_{s=t}^{\infty}(1+\gamma_{s})\right)V_{2}(t)+\sum_{s=t}^{\infty}\hat{\gamma}_{s},~~~~\forall t\geq t_{2}.
		\end{align}}
		\noindent By \eqref{eq:l2_pr_19}, it can be shown that $\{W(t)\}$ satisfies,
		{\small\begin{align}
		\label{eq:l2_pr_20}
		\mathbb{E}_{\theta^{*}}[W(t+1) | \mathcal{F}_{t}] \le W(t).
		\end{align}}
		\noindent Hence,  $\{W(t)\}$ is a non-negative super martingale and converges a.s. to a bounded random variable $W^{*}$ as $t\to\infty$. It then follows from \eqref{eq:l2_pr_19} that $V_{2}(t)\to W^{*}$ as $t\to\infty$. Thus, we conclude that the sequences $\{\mathbf{\theta}_{n}(t)\}$ are bounded for all $n$.
	\end{IEEEproof}
	\noindent Due to inherent stochasticity associated with the noisy observations, there need not be {\it uniform} boundedness of the estimate sequences.  Hence, while Lemma \ref{le:l0} establishes the pathwise boundedness of the parameter estimate sequence, it does not guarantee uniform boundedness over almost all sample paths.
	\begin{Lemma}
		\label{le:l1}
		Let the hypotheses of Theorem \ref{th:cons} hold. Then, we have,
		{\small\begin{align}
		\label{eq:l1_1}
		\mathbb{P}_{\btheta}\left(\lim_{t\to\infty}\mathbf{x}_{n}(t)=\btheta\right)=1,~\forall n.
		\end{align}}
	\end{Lemma}
	\begin{IEEEproof}	
		\noindent Denote the processes $\{\mathbf{z}_{n}(t)\}$ and $\{\widehat{\mathbf{z}}_{n}(t)\}$ as
		{\small\begin{align}
		\label{eq:l1_pr_2}
		\mathbf{z}_{n}(t)=\mathbf{x}_{n}(t)-\btheta~~\mbox{and}~~
		\widehat{\mathbf{z}}_{n}(t)=\wx_{n}(t)-\btheta
		\end{align}}
		respectively.
		\noindent Let {\small$\mathbf{z}(t)=\left[\mathbf{z}_{1}^{\top}(t)\cdots\mathbf{z}_{n}^{\top}(t)\right]^{\top}$} and {\small$\widehat{\mathbf{z}}(t)=\left[\widehat{\mathbf{z}}_{1}^{\top}(t)\cdots\widehat{\mathbf{z}}_{n}^{\top}(t)\right]^{\top}$}.
		\noindent From, \eqref{eq:dist4}, we have,
		{\small\begin{align}
		\label{eq:l1_pr_3}
		&\widehat{\mathbf{z}}(t+1)=\mathbf{z}(t)-\beta_{t}\left(\mathbf{L}\otimes\mathbf{I}_{M}\right)\mathbf{z}(t)\nonumber\\&+\alpha_{t}\mathbf{G}(\mathbf{x}(t))\mathbf{R}^{-1}\left(\mathbf{y}(t)-\mathbf{f}\left(\mathbf{x}(t)\right)\right),
		\end{align}}
		\noindent where we have used the fact that {\small$\left(\mathbf{L}\otimes\mathbf{I}_{M}\right)\left(\mathbf{1}_{N}\otimes\btheta\right)=\mathbf{0}$}.
		\noindent Define the $\{\mathcal{F}_{t}\}$-adapted process $\{V(t)\}$  by
		{\small\begin{align}
		\label{eq:l1_pr_4}
		V(t)=\left\|\mathbf{z}(t)\right\|^{2}.
		\end{align}}
		\noindent Now, using \eqref{eq:l0_pr_1} and the fact that {\small$\mathbb{E}_{\btheta}\left[\mathbf{y}(t)-\mathbf{f}\left(\mathbf{1}_{N}\otimes\btheta\right)\right]=\mathbf{0}$},
		\noindent we have,
		{\small\begin{align}
		\label{eq:l1_pr_5}
		&\mathbb{E}_{\btheta}[V(t+1)|\mathcal{F}_{t}]\le V(t) + \beta_{t}^{2}\mathbf{z}^{\top}(t)\left(\overline{\mathbf{L}} \otimes \mathbf{I}_{M}\right)^{2}\mathbf{z}(t)\nonumber\\
		&+\alpha_{t}^{2}\mathbb{E}_{\btheta}\left[\left\|\mathbf{G}\left(\mathbf{x}(t)\right)\mathbf{R}^{-1}\left(\mathbf{y}(t)-\mathbf{f}\left(\mathbf{1}_{N}\otimes\btheta\right)\right)\right\|^{2}\right]\nonumber\\&-2\beta_{t}\mathbf{z}^{\top}(t)\left(\overline{\mathbf{L}}\otimes \mathbf{I}_{M}\right)\mathbf{z}(t)\nonumber\\
		&-2\alpha_{t}\mathbf{z}(t)^{\top}\mathbf{G}\left(\mathbf{x}(t)\right)\mathbf{R}^{-1}\left(\mathbf{f}\left(\mathbf{x}(t)\right)-\mathbf{f}\left(\mathbf{1}_{N}\otimes\btheta\right)\right)\nonumber\\&+2\alpha_{t}\beta_{t}\mathbf{z}(t)^{\top}\left(\overline{\mathbf{L}} \otimes \mathbf{I}_{M}\right)\mathbf{G}\left(\mathbf{x}(t)\right)\mathbf{R}^{-1}\left(\mathbf{f}\left(\mathbf{x}(t)\right)-\mathbf{f}\left(\mathbf{1}_{N}\otimes\btheta\right)\right) \nonumber\\
		&+\alpha_{t}^{2}\left\|\left(\mathbf{f}\left(\mathbf{x}(t)\right)-\mathbf{f}\left(\mathbf{1}_{N}\otimes\btheta\right)\right)^{\top}\mathbf{G}\left(\mathbf{x}(t)\right)^{\top}\mathbf{R}^{-1}\right\|^{2}.
		\end{align}}
		\noindent Following the steps as in the proof of Lemma \ref{l2} and using \eqref{eq:l0_pr_9}-\eqref{eq:l0_pr_14}, we have,
		{\small\begin{align}
		\label{eq:l1_pr_13}
		&\mathbb{E}_{\btheta}[V(t+1)|\mathcal{F}_{t}]\le\left(1-c_{7}\alpha_{t}\right)V(t)+c_{7}\alpha_{t}^{2}\nonumber\\
		&\le V(t)+c_{7}\alpha_{t}^{2},
		\end{align}}
		\noindent for an appropriately chosen positive constant $c_{7}$.
		\noindent Now, consider the $\{\mathcal{F}_{t}\}$-adapted process $\{V_{1}(t)\}$ defined as follows
		{\small\begin{align}
		\label{eq:l1_pr_14}
		&V_{1}(t)=V(t)-c_{7}\sum_{s=t}^{\infty}\alpha_{s}^{2}\nonumber\\
		&=V(t)-bc_{7}\sum_{s=t}^{\infty}(t+1)^{-2}.
		\end{align}}
		\noindent Since, $\{(t+1)^{-2}\}$ is summable, the process $\{V_{1}(t)\}$ is bounded from below. Moreover, it also follows that $\{V_{1}(t)\}_{t\geq t_{1}}$ is a supermartingale and hence converges a.s. to a finite random variable. By definition from \eqref{eq:l1_pr_14}, we also have that $\{V(t)\}$ converges to a non-negative finite random variable $V^{*}$. Finally, from \eqref{eq:l1_pr_13}, we have that,
		{\small\begin{align}
		\label{eq:l1_pr_15}
		\mathbb{E}_{\btheta}[V(t+1)]\le \left(1-c_{1}\alpha_{t}\right)\mathbb{E}_{\btheta}[V(t)]+bc_{7}(t+1)^{-2},
		\end{align}}
		\noindent for $t\geq t_{1}$. The sequence $\{V(t)\}$ then falls under the purview of Lemmas 4 and 5 of \cite{kar2011convergence}, and we have $\mathbb{E}_{\btheta}[V(t)]\to 0$ as $t\to\infty$. Finally, by Fatou's Lemma, where we use the non-negativity of the sequence $\{V(t)\}$, we conclude that
		{\small\begin{align}
		\label{eq:l1_pr_16}
		0\leq \mathbb{E}_{\btheta}[V^{*}]\le\liminf_{t\to\infty}\mathbb{E}_{\btheta}[V(t)]=0,
		\end{align}}
		\noindent which thus implies that $V^{*}=0$ a.s. Hence, $\left\|\mathbf{z}(t)\right\|\to 0$ as $t\to\infty$ and the desired assertion follows.
	\end{IEEEproof}	
	\noindent We will use the following approximation result~(Lemma \ref{l33}) and the generalized convergence criterion~(Lemma \ref{l34}) for the proof of Theorem \ref{th:cons}.
	\begin{Lemma}[Lemma 4.3 in \cite{Fabian-1}]
		\label{l33}
		Let $\{b_{t}\}$ be a scalar sequence satisfying
		{\small\begin{align}
		\label{eq:l33_1}
		b_{t+1}\le \left(1-\frac{c}{t+1}\right)b_{t}+d_{t}(t+1)^{-\tau},
		\end{align}}
		where $c > \tau, \tau > 0$, and the sequence $d_{t}$ is summable. Then, we have,
		{\small\begin{align}
		\label{eq:l33_2}
		\limsup_{t\to\infty}~(t+1)^{\tau}b_{t}<\infty.
		\end{align}}
	\end{Lemma}
	\begin{Lemma}[Lemma 10 in \cite{dubins1965sharper}]
		\label{l34}
		Let $\{J(t)\}$ be an $\mathbb{R}$-valued $\{\mathcal{F}_{t+1}\}$-adapted process such that $\mathbb{E}\left[J(t)|\mathcal{F}_{t}\right]=0$ a.s. for each $t\geq 1$. Then the sum $\sum_{t\geq 0}J(t)$ exists and is finite a.s. on the set where $\sum_{t\geq 0}\mathbb{E}\left[J(t)^{2}|\mathcal{F}_{t}\right]$ is finite.
	\end{Lemma}
	\noindent We now return to the proof of Theorem \ref{th:cons}.
	\noindent Define $\bar{\tau}\in[0,1/2)$ such that,
	{\small\begin{align}
	\label{eq:t1_pr_1}
	\mathbb{P}_{\btheta}\left(\lim_{t\to\infty}(t+1)^{\bar{\tau}}\left\|\mathbf{z}(t)\right\|=0\right)=1,
	\end{align}}
	\noindent where $\{\mathbf{z}(t)\}$ is as defined in \eqref{eq:l1_pr_2} and note that such a $\bar{\tau}$ always exists by Lemma \ref{le:l1}~(in particular $\bar{\tau}=0$). We now analyze and show that there exists a $\tau$ such that $\bar{\tau} < \tau <1/2$ for which the claim in \eqref{eq:t1_pr_1} holds. Now, choose a $\widehat{\tau} \in (\tau, 1/2)$ and let $\mu=(\widehat{\tau}+\bar{\tau})/2$. The recursion for $\{\mathbf{z}(t)\}$ can be written as follows:
	{\small\begin{align}
	\label{eq:t1_pr_2}
		&\left\|\mathbf{z}(t+1)\right\|^{2}\le\left\|\mathbf{z}(t)\right\|^{2}-2\beta_{t}\mathbf{z}^{\top}(t)\left(\mathbf{L} \otimes \mathbf{I}_{M}\right)\mathbf{z}(t) \nonumber\\&-2\alpha_{t}\mathbf{z}^{\top}(t)\mathbf{G}\left(\mathbf{x}(t)\right)\mathbf{R}^{-1}\left(\mathbf{f}\left(\mathbf{x}(t)\right)-\mathbf{f}\left(\mathbf{1}_{N}\otimes\btheta\right)\right)\nonumber\\
		&+\beta_{t}^{2}\mathbf{z}^{\top}(t)\left(\mathbf{L}\otimes \mathbf{I}_{M}\right)^{2}\mathbf{z}(t)\nonumber\\&+2\alpha_{t}\beta_{t}\mathbf{z}^{\top}(t)\left(\mathbf{L}\otimes \mathbf{I}_{M}\right)\mathbf{G}\left(\mathbf{x}(t)\right)\mathbf{R}^{-1}\left(\mathbf{f}\left(\mathbf{x}(t)\right)-\mathbf{f}\left(\mathbf{1}_{N}\otimes\theta\right)\right)\nonumber\\
		&+\alpha_{t}^{2}\left\|\mathbf{G}\left(\mathbf{x}(t)\right)\mathbf{R}^{-1}\left(\mathbf{y}(t)-\mathbf{f}\left(\mathbf{1}_{N}\otimes\btheta\right)\right)\right\|^{2}\nonumber\\&+\alpha_{t}^{2}\left\|\mathbf{G}\left(\mathbf{x}(t)\right)\mathbf{R}^{-1}\left(\mathbf{f}\left(\mathbf{x}(t)\right)-\mathbf{f}\left(\mathbf{1}_{N}\otimes\btheta\right)\right)\right\|^{2}\nonumber\\
		&+2\alpha_{t}\mathbf{z}^{\top}(t)\mathbf{G}\left(\mathbf{x}(t)\right)\mathbf{R}^{-1}\left(\mathbf{y}(t)-\mathbf{f}(\mathbf{1}_{N}\otimes\btheta)\right)\nonumber\\&+2\alpha_{t}^{2}\left(\mathbf{y}(t)-\mathbf{f}(\mathbf{1}_{N}\otimes\btheta)\right)^{\top}\mathbf{R}^{-1}\mathbf{G}^{\top}\left(\mathbf{x}(t)\right)\nonumber\\&\times\mathbf{G}\left(\mathbf{x}(t)\right)\mathbf{R}^{-1}\left(\mathbf{f}\left(\mathbf{1}_{N}\otimes\btheta\right)-\mathbf{f}\left(\mathbf{x}(t)\right)\right).
	\end{align}}
	
	\noindent Let {\small$\mathbf{J}(t)=\mathbf{G}\left(\mathbf{x}(t)\right)\mathbf{R}^{-1}\left(\mathbf{y}(t)-\mathbf{f}\left(\mathbf{1}_{N}\otimes\btheta\right)\right)$}.
	\noindent Now, we consider the term $\alpha_{t}^{2}\left\|\mathbf{J}(t)\right\|^{2}$. Since, the noise process under consideration has finite second moment and $2\mu < 1$, we have,
	
	{\small\begin{align}
	\label{eq:t1_pr_3}
	\sum_{t\geq 0}(t+1)^{2\mu}\alpha_{t}^{2}\left\|\mathbf{J}(t)\right\|^{2} < \infty.
	\end{align}}
	\noindent Let {\small$\mathbf{W}(t)=\alpha_{t}\mathbf{z}(t)^{\top}\mathbf{G}\left(\mathbf{x}(t)\right)\mathbf{R}^{-1}\left(\mathbf{y}(t)-\mathbf{f}(\mathbf{1}_{N}\otimes\btheta)\right)$}. It follows that $\mathbb{E}_{\btheta}\left[\mathbf{W}(t)|\mathcal{F}_{t}\right]=0$. We also have that\\ {\small$\mathbb{E}_{\btheta}\left[\mathbf{W}^{2}(t)|\mathcal{F}_{t}\right]\le\alpha_{t}^{2}\left\|\mathbf{z}(t)\right\|^{2}\left\|\mathbf{J}(t)\right\|^{2}$}. Noting that the noise under consideration has finite second order moment, we have,
	{\small\begin{align}
	\label{eq:t1_pr_4}
	\mathbb{E}_{\btheta}\left[\mathbf{W}^{2}(t)|\mathcal{F}_{t}\right]=o\left((t+1)^{-2-2\bar{\tau}}\right),
	\end{align}} and, hence,
	{\small\begin{align}
	\label{eq:t1_pr_5}
	\mathbb{E}_{\btheta}\left[(t+1)^{4\mu}\mathbf{W}^{2}(t)|\mathcal{F}_{t}\right]=o\left((t+1)^{-2+2\widehat{\tau}}\right).
	\end{align}}
	\noindent Hence, by Lemma \ref{l34}, we conclude that {\small$\sum_{t\geq 0}(t+1)^{2\mu}\mathbf{W}(t)$} exists and is finite, as $2\widehat{\tau} < 1$. Similarly, it can be shown that, for {\small$\mathbf{W}_{1}(t)=\alpha_{t}\beta_{t}\mathbf{z}(t)^{\top}\mathbf{G}\left(\mathbf{x}(t)\right)\mathbf{R}^{-1}\left(\mathbf{f}(\mathbf{1}_{N}\otimes\btheta)-\mathbf{y}(t)\right)$}, the sum $\sum_{t\geq 0}(t+1)^{2\mu}\mathbf{W}_{1}(t)$ exists and is finite. Finally, consider {\small$\mathbf{W}_{2}(t)=\alpha_{t}^{2}\left(\mathbf{y}(t)-\mathbf{f}\left(\mathbf{1}_{N}\otimes\btheta\right)\right)^{\top}\mathbf{R}^{-1}\mathbf{G}\left(\mathbf{x}(t)\right)^{\top}\times\mathbf{G}\left(\mathbf{x}(t)\right)\mathbf{R}^{-1}\left(\mathbf{f}\left(\mathbf{1}_{N}\otimes\btheta\right)-\mathbf{f}\left(\mathbf{x}(t)\right)\right)$}. It follows that {\small$\mathbb{E}_{\btheta}\left[\mathbf{W}_{2}(t)|\mathcal{F}_{t}\right]=0$}. We also have that {\small$\mathbb{E}_{\btheta}\left[\mathbf{W}_{2}^{2}(t)|\mathcal{F}_{t}\right]\le\alpha_{t}^{4}\left\|\mathbf{z}(t)\right\|^{2}\left\|\mathbf{J}(t)\right\|^{2}$}. Following as in \eqref{eq:t1_pr_4} and \eqref{eq:t1_pr_5}, we have that $\sum_{t\geq 0}(t+1)^{2\mu}\mathbf{W}_{2}(t)$ exists and is finite. Using all the inequalities derived in \eqref{eq:l0_pr_9}-\eqref{eq:l0_pr_14}, we have,
	{\small\begin{align}
	\label{eq:t1_pr_6}
	&\left\|\mathbf{z}(t+1)\right\|^{2}\le\left(1-c_{1}\alpha_{t}+c_{2}\alpha_{t}\beta_{t}+c_{3}\alpha_{t}^{2}\right)\left\|\mathbf{z}(t)\right\|^{2}+\alpha_{t}^{2}\left\|\mathbf{J}(t)\right\|^{2}\nonumber\\&-c_{6}(\beta_{t}-\beta^{2}_{t})\left\|\mathbf{z}_{C\perp}(t)\right\|^{2}+2\mathbf{W}(t)+2\mathbf{W}_{1}(t)+2\mathbf{W}_{2}(t).
	\end{align}}
	\noindent Finally, noting that $c_{1}\alpha_{t}$ dominates $c_{2}\alpha_{t}\beta_{t}$ and $c_{3}\alpha_{t}^{2}$, $\beta_{t}$ dominates $\beta_{t}^{2}$, we have eventually
	{\small\begin{align}
	\label{eq:t1_pr_7}
	&\left\|\mathbf{z}(t+1)\right\|^{2}\le\left(1-c_{1}\alpha_{t}\right)\left\|\mathbf{z}(t)\right\|^{2}\nonumber\\&+\alpha_{t}^{2}\left\|\mathbf{J}(t)\right\|^{2}+2\mathbf{W}(t)+2\mathbf{W}_{1}(t)+2\mathbf{W}_{2}(t).
	\end{align}}
	\noindent To this end, using the analysis in \eqref{eq:t1_pr_3}-\eqref{eq:t1_pr_5}, we have, from \eqref{eq:t1_pr_7}
	{\small\begin{align}
	\label{eq:t1_pr_8}
	\left\|\mathbf{z}(t+1)\right\|^{2}\le\left(1-c_{1}\alpha_{t}\right)\left\|\mathbf{z}(t)\right\|^{2}+d_{t}(t+1)^{-2\mu},
	\end{align}}
	\noindent where
	{\small\begin{align}
	\label{eq:t1_pr_9}
	d_{t}(t+1)^{-2\mu}=\alpha_{t}^{2}\left\|\mathbf{J}(t)\right\|^{2}+2\mathbf{W}(t)++2\mathbf{W}_{1}(t)+2\mathbf{W}_{2}(t).
	\end{align}}
	\noindent Finally, noting that $c_{1}\alpha_{t}(t+1)=1>2\mu$, an immediate application of Lemma \ref{l33} gives
	{\small\begin{align}
	\label{eq:t1_pr_10}
	\limsup_{t\to\infty}(t+1)^{2\mu}\left\|\mathbf{z}(t)\right\|^{2} < \infty~ a.s.
	\end{align}}
	\noindent So, we have that there exists a $\tau$ with $\bar{\tau}<\tau<\mu$ for which $(t+1)^{\tau}\left\|\mathbf{z}(t)\right\|\to 0$ as $t\to\infty$. Thus, for every $\bar{\tau}$ for which \eqref{eq:th_cons1} holds, there exists $\tau\in(\bar{\tau},1/2)$ for which the result in \eqref{eq:th_cons1} still continues to hold. By a simple application of induction, we conclude that the result holds for all $\tau \in [0,1/2)$.	
\end{IEEEproof}

\subsection{Proof of Theorem \ref{th:2}}
\label{subsec:th_2}
\begin{IEEEproof}
	\noindent The proof of Theorem \ref{th:2} needs the following Lemma from \cite{Fabian-2} concerning the asymptotic normality of the stochastic recursions.
	\begin{Lemma}[Theorem 2.2 in \cite{Fabian-2}]
		\label{main_res_l0}
		\noindent Let $\{\mathbf{z}_{t}\}$ be an $\mathbb{R}^{k}$-valued $\{\mathcal{F}_{t}\}$-adapted process that satisfies
		{\small\begin{align}
		\label{eq:l0_1}
		\mathbf{z}_{t+1}=\left(\mathbf{I}_{k}-\frac{1}{t+1}\Gamma_{t}\right)\mathbf{z}_{t}+(t+1)^{-1}\mathbf{\Phi}_{t}\mathbf{V}_{t}+(t+1)^{-3/2}\mathbf{T}_{t},
		\end{align}}
		\noindent where the stochastic processes {\small$\{\mathbf{V}_{t}\}, \{\mathbf{T}_{t}\} \in \mathbb{R}^{k}$} while {\small$\{\mathbf{\Gamma}_{t}\}, \{\mathbf{\Phi}_{t}\} \in \mathbb{R}^{k\times k}$}. Moreover, suppose for each $t$, $\mathbf{V}_{t-1}$ and $\mathbf{T}_{t}$ are $\mathcal{F}_{t}$-adapted, whereas the processes {\small$\{\mathbf{\Gamma}_{t}\}$, $\{\mathbf{\Phi}_{t}\}$ are $\{\mathcal{F}_{t}\}$}-adapted.
		
		\noindent Also, assume that
		{\small\begin{align}
		\label{eq:l0_2}
		\mathbf{\Gamma}_{t}\to\mathbf{\Gamma}, \mathbf{\Phi}_{t}\to\mathbf{\Phi},~ \textit{and} ~\mathbf{T}_{t}\to 0 ~~\mbox{a.s. as $t\rightarrow\infty$},
		\end{align}}
		\noindent where {\small$\mathbf{\Gamma}$} is a symmetric and positive definite matrix, and admits an eigen decomposition of the form {\small$\mathbf{P}^{\top}\mathbf{\Gamma}\mathbf{P}=\mathbf{\Lambda}$}, where {\small$\mathbf{\Lambda}$} is a diagonal matrix and {\small$\mathbf{P}$} is an orthogonal matrix. Furthermore, let the sequence {\small$\{\mathbf{V}_{t}\}$} satisfy {\small$\mathbb{E}\left[\mathbf{V}_{t}|\mathcal{F}_{t}\right]=0$} for each $t$ and suppose there exists a positive constant $C$ and a matrix $\Sigma$ such that {\small$C > \left\|\mathbb{E}\left[\mathbf{V}_{t}\mathbf{V}_{t}^{\top}|\mathcal{F}_{t}\right]-\Sigma\right\|\to 0~a.s. ~\textit{as}~ t\to\infty$} and with {\small$\sigma_{t,r}^{2}=\int_{\left\|\mathbf{V}_{t}\right\|^{2} \ge r(t+1)}\left\|\mathbf{V}_{t}\right\|^{2}d\mathbb{P}$, let $\lim_{t\to\infty}\frac{1}{t+1}\sum_{s=0}^{t}\sigma_{s,r}^{2}=0$} for every $r > 0$. Then, we have,
		{\small\begin{align}
		\label{eq:l0_3}
		(t+1)^{1/2}\mathbf{z}_{t}\overset{\mathcal{D}}{\Longrightarrow}\mathcal{N}\left(\mathbf{0}, \mathbf{P}\mathbf{M}\mathbf{P}^{\top}\right),
		\end{align}}
		\noindent where the $(i,j)$-th entry of the matrix $\mathbf{M}$ is given by
		{\small\begin{align}
		\label{eq:l0_4}
		\left[\mathbf{M}\right]_{ij}=\left[\mathbf{P}^{\top}\mathbf{\Phi}\mathbf{\Sigma}\mathbf{\Phi}^{\top}\mathbf{P}\right]_{ij}\left(\left[\mathbf{\Lambda}\right]_{ii}+\left[\mathbf{\Lambda}\right]_{jj}-1\right)^{-1}.
		\end{align}}
	\end{Lemma}
	\noindent From Theorem \ref{th:cons} and the fact that $\btheta$ lies in the interior of the parameter set $\Theta$, we have that there exists an $\epsilon > 0$, such that $B_{\epsilon}(\btheta) \in \Theta$, where $B_{\epsilon}(\btheta)$ denotes the open ball centered at $\btheta$ with radius $\epsilon$. In particular, fix an $\epsilon > 0$ for which $B_{\epsilon}(\btheta) \in \Theta$. Then, we have that there exists a random time $t_{\epsilon}(\omega)$, which is almost surely finite, i.e., $\mathbb{P}\left(t_{\epsilon}(\omega)<\infty\right)=1$, such that $\left\|\mathbf{x}_{\omega}(t)-\btheta\right\| < \epsilon$ for all $t\geq t_{\epsilon}(\omega)$, where $\omega$ denotes the sample path. In the above, we introduce the $\omega$-argument to emphasize that the time $t_{\epsilon}(\omega)$ is random and sample-path dependent and our analysis is pathwise. With the above development in place, we note that \eqref{eq:dist1} and \eqref{eq:dist2} can be rewritten as follows:
	{\small\begin{align}
	\label{eq:t2_pr_0}
	&\mathbf{x}(t+1)=\mathbf{x}(t)-\beta_{t}\left(\mathbf{L}\otimes\mathbf{I}_{M}\right)\mathbf{x}(t)\nonumber\\&+\alpha_{t}\mathbf{G}(\mathbf{x}(t))\mathbf{R}^{-1}\left(\mathbf{y}(t)-\mathbf{f}\left(\mathbf{x}(t)\right)\right)+\mathbf{e}_{\mathcal{P}}(t),~\forall t\geq 0,
	\end{align}}
	\noindent where $\mathbf{e}_{\mathcal{P}}(t)$ is the projection error which is given by,
	{\small\begin{align}
	\label{eq:t2_pr_1}
	\mathbf{e}_{\mathcal{P}}(t)=\mathbf{x}(t+1)-\wx(t+1),
	\end{align}}
	\noindent and in particular $\mathbf{e}_{\mathcal{P}}(t)=0$ for all $t\ge t_{\epsilon}$.

\noindent Define the process $\left\{\mathbf{x}_{\mbox{\scriptsize{avg}}}(t)\right\}$ as
{\small\begin{align}
\label{eq:t2_pr_2}
\mathbf{x}_{\mbox{\scriptsize{avg}}}(t)=\left(\frac{\mathbf{1}_{N}^{\top}}{N}\otimes\mathbf{I}_{M}\right)\mathbf{x}(t).
\end{align}}
\noindent It is readily seen that the process $\left\{\mathbf{x}_{\mbox{\scriptsize{avg}}}(t)\right\}$ satisfies the recursion
\vspace{-10pt}
{\small\begin{align}
	\label{eq:t2_pr_3}
	&\mathbf{x}_{\mbox{\scriptsize{avg}}}(t+1)=\mathbf{x}_{\mbox{\scriptsize{avg}}}(t)+\left(\frac{\mathbf{1}_{N}^{\top}}{N}\otimes\mathbf{I}_{M}\right)\mathbf{e}_{\mathcal{P}}(t)\nonumber\\&+\alpha_{t}\left(\frac{\mathbf{1}_{N}^{\top}}{N}\otimes\mathbf{I}_{M}\right)\mathbf{G}(\mathbf{x}(t))\mathbf{R}^{-1}\left(\mathbf{y}(t)-\mathbf{f}\left(\mathbf{x}(t)\right)\right)\nonumber\\
	&=\mathbf{x}_{\mbox{\scriptsize{avg}}}(t)+\left(\frac{\mathbf{1}_{N}^{\top}}{N}\otimes\mathbf{I}_{M}\right)\mathbf{e}_{\mathcal{P}}(t)\nonumber\\&+\frac{\alpha_{t}}{N}\sum_{n=1}^{N}\nabla \mathbf{f}_{n}\left(\mathbf{x}_{n}(t)\right)\mathbf{R}_{n}^{-1}\left(\mathbf{y}_{n}(t)-\mathbf{f}_{n}\left(\mathbf{x}_{n}(t)\right)\right).
	\end{align}}
\noindent Now noting that, for all $t\geq 0$, $\mathbf{x}_{n}(t)\in\Theta$ for each $n$ and as $\Theta$ is a convex set, we have that $\mathbf{x}_{\mbox{\scriptsize{avg}}}(t)\in\Theta$ for all $t\geq 0$. Then, we have by the mean-value theorem for each agent $n$
\vspace{-10pt}
{\small\begin{align}
\label{eq:t2_pr_4}
&\mathbf{f}_{n}\left(\mathbf{x}_{\mbox{\scriptsize{avg}}}(t)\right)=\mathbf{f}_{n}\left(\btheta\right)+\nabla^{\top}\mathbf{f}_{n}\left(c\btheta+(1-c)\mathbf{x}_{\mbox{\scriptsize{avg}}}(t)\right)\left(\mathbf{x}_{\mbox{\scriptsize{avg}}}(t)-\btheta\right),
\end{align}}
\noindent where $0<c<1$. It is to be noted that {\small$\nabla^{\top}\mathbf{f}_{n}\left(c\btheta+(1-c)\mathbf{x}_{\mbox{\scriptsize{avg}}}(t)\right)\to\nabla^{\top}\mathbf{f}_{n}\left(\btheta\right)$} as {\small$\mathbf{x}_{\mbox{\scriptsize{avg}}}(t)\to\btheta$} in the limit $t\to\infty$.
\noindent Using \eqref{eq:t2_pr_4} in \eqref{eq:t2_pr_3}, we have for all $t\geq 0$,
{\small\begin{align}
	\label{eq:t2_pr_5}
	&\mathbf{x}_{\mbox{\scriptsize{avg}}}(t+1)-\btheta=\left(\mathbf{I}-\frac{\alpha_{t}}{N}\left(\sum_{n=1}^{N}\nabla\mathbf{f}_{n}\left(\mathbf{x}_{n}(t)\right)\mathbf{R}_{n}^{-1}\right.\right.\nonumber\\&\times\left.\left.\nabla^{\top}\mathbf{f}_{n}\left(c\btheta+(1-c)\mathbf{x}_{\mbox{\scriptsize{avg}}}(t)\right)\right)\right)\left(\mathbf{x}_{\mbox{\scriptsize{avg}}}(t)-\btheta\right)\nonumber\\&+\frac{\alpha_{t}}{N}\sum_{n=1}^{N}\nabla\mathbf{f}_{n}\left(\mathbf{x}_{n}(t)\right)\mathbf{R}_{n}^{-1}\mathbf{\zeta}_{n}(t)+\left(\frac{\mathbf{1}_{N}^{\top}}{N}\otimes\mathbf{I}_{M}\right)\mathbf{e}_{\mathcal{P}}(t)\nonumber\\&+\frac{\alpha_{t}}{N}\sum_{n=1}^{N}\nabla\mathbf{f}_{n}\left(\mathbf{x}_{n}(t)\right)\mathbf{R}_{n}^{-1}\left(\mathbf{f}_{n}\left(\mathbf{x}_{\mbox{\scriptsize{avg}}}(t)\right)-\mathbf{f}_{n}\left(\mathbf{x}_{n}(t)\right)\right).
	\end{align}}
\noindent The following Lemma will be crucial for the subsequent part of the proof.
\begin{Lemma}
	\label{l_t2_pr}
	For every $\tau_{0}$ such that $0\le\tau_{0} < 1-\tau_{1}-1/(2+\epsilon_{1})$, we have,
	\vspace{-10pt}
	{\small\begin{align}
	\mathbb{P}_{\btheta}\left(\lim_{t\to\infty}(t+1)^{\tau_{0}}\left(\mathbf{x}_{n}(t)-\mathbf{x}_{\mbox{\scriptsize{avg}}}(t)\right)=0\right)=1.
	\end{align}}
\end{Lemma}
\begin{IEEEproof}
	\noindent The proof follows exactly like the proof of Lemma IV.2 in \cite{kar2014asymptotically}. Note the additional term that comes up in $\mathbf{x}_{n}(t)-\mathbf{x}_{\mbox{\scriptsize{avg}}}(t)$ in the current context due to the projection error is given by {\small$\mathbf{e}_{\mathcal{P},n}(t)-\left(\frac{\mathbf{1}_{N}^{\top}}{N}\otimes\mathbf{I}_{M}\right)\mathbf{e}_{\mathcal{P}}(t)$}; nonetheless, this term satisfies the property that {\small$\mathbb{P}_{\btheta}\left(\lim_{t\to\infty}(t+1)^{\tau_{0}}\left(\mathbf{e}_{\mathcal{P},n}(t)-\left(\frac{\mathbf{1}_{N}^{\top}}{N}\otimes\mathbf{I}_{M}\right)\mathbf{e}_{\mathcal{P}}(t)\right)=0\right)=1$} as {\small$\mathbf{e}_{\mathcal{P}}(t)=0$} for all $t\ge t_{\epsilon}$. Hence, the techniques employed in the proof of Lemma IV.2 also apply here. Lemma IV.2 in \cite{kar2014asymptotically} is concerned with the asymptotic agreement of the estimates across any pair of agents, but as an intermediate result the agreement of the estimate at an agent and the averaged estimate is established.
\end{IEEEproof}
	\noindent As $\tau_{1}+1/(2+\epsilon_{1}) < 1/2$, from Lemma \ref{l_t2_pr} we have that there exists an $\epsilon_{2} > 0$ (sufficiently small) such that
	{\small\begin{align}
	\label{eq:t2_pr_6}
	\mathbb{P}_{\btheta}\left(\lim_{t\to\infty}(t+1)^{\frac{1}{2}+\epsilon_{2}}\left(\mathbf{x}_{n}(t)-\mathbf{x}_{\mbox{\scriptsize{avg}}}(t)\right)=0\right)=1.
	\end{align}}
	\noindent We consider the process $\left\{\mathbf{x}_{\mbox{\scriptsize{avg}}}(t)\right\}$ for the application of Lemma \ref{main_res_l0}.
	\noindent Hence, comparing term by term of \eqref{eq:t2_pr_5} with \eqref{eq:l0_1}, we have,
	{\small\begin{align}
		\label{eq:t2_pr_7}
		&\mathbf{\Gamma}_{t}=\frac{a}{N}\sum_{n=1}^{N}\nabla\mathbf{f}_{n}\left(\mathbf{x}_{n}(t)\right)\mathbf{R}_{n}^{-1}\nabla^{\top}\mathbf{f}_{n}\left(c\btheta+(1-c)\mathbf{x}_{\mbox{\scriptsize{avg}}}(t)\right)\nonumber\\&\to \frac{a}{N}\sum_{n=1}^{N}\nabla\mathbf{f}_{n}\left(\btheta\right)\mathbf{R}_{n}^{-1}\nabla\mathbf{f}_{n}^{\top}\left(\btheta\right)=a\mathbf{\Gamma}_{\btheta},\nonumber\\
		&\mathbf{\Phi}_{t}=a\left(\frac{\mathbf{1}_{N}^{\top}}{N}\otimes\mathbf{I}_{M}\right)\mathbf{G}(\mathbf{x}(t))\mathbf{R}^{-1}\to a\left(\frac{\mathbf{1}_{N}^{\top}}{N}\otimes\mathbf{I}_{M}\right)\mathbf{G}(\mathbf{1}\otimes\btheta)\mathbf{R}^{-1}\nonumber\\&=a\mathbf{\Phi},\nonumber\\
		&\mathbf{V}_{t}=\mathbf{\zeta}(t), \mathbb{E}\left[\mathbf{V}_{t}|\mathcal{F}_{t}\right]=0, \mathbb{E}\left[\mathbf{V}_{t}\mathbf{V}^{\top}_{t}|\mathcal{F}_{t}\right]=\mathbf{R},\nonumber\\
		&\mathbf{T}_{t}=a(t+1)^{1/2}\left(\frac{\mathbf{1}_{N}^{\top}}{N}\otimes\mathbf{I}_{M}\right)\mathbf{G}(\mathbf{x}(t))\mathbf{R}^{-1}\times\nonumber\\&\left(\mathbf{f}\left(\mathbf{1}\otimes\mathbf{x}_{\mbox{\scriptsize{avg}}}(t)\right)-\mathbf{f}\left(\mathbf{x}(t)\right)\right)+(t+1)^{3/2}\left(\frac{\mathbf{1}_{N}^{\top}}{N}\otimes\mathbf{I}_{M}\right)\mathbf{e}_{\mathcal{P}}(t)\to 0,
		\end{align}}
	\noindent where the convergence of $\mathbf{T}_{t}$ follows from Lemma \ref{l_t2_pr}. Due to the i.i.d nature of the noise process, we have the uniform integrability condition for the process $\{\mathbf{V}_{t}\}$. Hence, $\{\mathbf{x}_{\mbox{\scriptsize{avg}}}(t)\}$ falls under the purview of Lemma \ref{main_res_l0}, and we thus conclude that
	{\small\begin{align}
	\label{eq:t2_pr_8}
	(t+1)^{1/2}\left(\mathbf{x}_{\mbox{\scriptsize{avg}}}(t)-\btheta\right)\overset{\mathcal{D}}{\Longrightarrow}\mathcal{N}(0,\mathbf{P}\mathbf{M}\mathbf{P}^{\top}),
	\end{align}}
	\noindent where
	{\small\begin{align}
	\label{eq:t2_pr_9}
	&a\mathbf{P}^{\top}\mathbf{\Gamma}_{\btheta}\mathbf{P}=a\mathbf{\Lambda}_{\btheta},\nonumber\\
	&\left[\mathbf{M}\right]_{ij}=\left[a^{2}\mathbf{P}^{\top}\mathbf{\Phi}\mathbf{R}\mathbf{\Phi}^{\top}\mathbf{P}\right]_{ij}\left(a\left[\mathbf{\Lambda}\right]_{\btheta,ii}+a\left[\mathbf{\Lambda}\right]_{\btheta,jj}-1\right)^{-1}\nonumber\\&=\frac{a^{2}}{N}\left[\mathbf{\Lambda}\right]_{ij}\left(a\left[\mathbf{\Lambda}\right]_{\btheta,ii}+a\left[\mathbf{\Lambda}\right]_{\btheta,jj}-1\right)^{-1},
	\end{align}}
	\noindent which also implies that $\mathbf{M}$ is a diagonal matrix with its $i$-th diagonal element given by {\small$\frac{a^{2}\mathbf{\Lambda}_{\btheta,ii}}{2aN\mathbf{\Lambda}_{\btheta,ii}-N}$}. We already have that {\small$\mathbf{P}\mathbf{\Lambda}_{\btheta}\mathbf{P}^{\top}=\mathbf{\Gamma}_{\btheta}$}. Hence, the matrix with eigenvalues as {\small$\frac{a^{2}\mathbf{\Lambda}_{\btheta,ii}}{2aN\mathbf{\Lambda}_{\btheta,ii}-N}$} is given by
	{\small\begin{align}
	\label{eq:t2_pr_10}
	\mathbf{P}\mathbf{M}\mathbf{P}^{\top}=\frac{a\mathbf{I}}{2N}+\frac{\left(\mathbf{\Gamma}_{\btheta}-\frac{\mathbf{I}}{2a}\right)^{-1}}{4N}.
	\end{align}}
	\noindent Now from \eqref{eq:t2_pr_6}, which is a consequence of Lemma \ref{l_t2_pr}, we have that the processes $\{\mathbf{x}_{n}(t)\}$ and $\{\mathbf{x}_{\mbox{\scriptsize{avg}}}(t)\}$ are indistinguishable in the $t^{1/2}$ time scale, which is formalized as follows:
	{\small\begin{align}
	\label{eq:t2_pr_11}
	&\mathbb{P}_{\btheta}\left(\lim_{t\to\infty}\left\|\sqrt{t+1}\left(\mathbf{x}_{n}(t)-\btheta\right)-\sqrt{t+1}\left(\mathbf{x}_{\mbox{\scriptsize{avg}}}(t)-\btheta\right)\right\|=0\right)\nonumber\\
	&=\mathbb{P}_{\btheta}\left(\lim_{t\to\infty}\left\|\sqrt{t+1}\left(\mathbf{x}_{n}(t)-\mathbf{x}_{\mbox{\scriptsize{avg}}}(t)\right)\right\|=0\right)=1.
	\end{align}}
	\noindent Thus, the difference of the sequences {\small$\left\{\sqrt{t+1}\left(\mathbf{x}_{n}(t)-\btheta\right)\right\}$ and $\left\{\sqrt{t+1}\left(\mathbf{x}_{\mbox{\scriptsize{avg}}}(t)-\btheta\right)\right\}$} converges a.s. to zero as $t\rightarrow\infty$, and hence we have,
	{\small\begin{align}
	\label{eq:t2_pr_12}
	\sqrt{t+1}\left(\mathbf{x}_{n}(t)-\btheta\right)\overset{\mathcal{D}}{\Longrightarrow}\mathcal{N}\left(0,\frac{a\mathbf{I}}{2N}+\frac{\left(\mathbf{\Gamma}_{\btheta}-\frac{\mathbf{I}}{2a}\right)^{-1}}{4N}\right).
	\end{align}}
\end{IEEEproof}
\vspace{-10pt}
\section{Conclusion}
\label{conclusion} In this paper, we have considered the problem of distributed recursive parameter estimation in a network of sparsely interconnected agents. We have proposed a \emph{consensus} + \emph{innovations} nonlinear least squares type algorithm, $\mathcal{CIWNLS}$, in which every agent updates its parameter estimate at every observation sampling epoch by simultaneous processing of neighborhood information and locally sensed new information and in which the inter-agent collaboration is restricted to a possibly sparse but connected communication graph. Under rather weak conditions, connectivity of the inter-agent communication and a \emph{global observability} criterion, we have shown that the proposed algorithm leads to consistent  parameter estimates at each agent. Furthermore, under standard smoothness assumptions on the sensing nonlinearities, we have established order-optimal pathwise convergence rates and the asymptotic normality of the parameter estimate sequences generated by the proposed distributed estimator $\mathcal{CIWNLS}$. A natural direction for future research consists of obtaining techniques and conditions to obtain innovation gains so as to reduce the gap between the agent asymptotic covariances and that of the centralized WNLS estimator. Methods developed in~\cite{kar2014asymptotically} may be employed and extended to obtain such characterization. The adaptation of specialized distributed algorithms, for example, distributed localization studied in~\cite{khan2009distributed} to the framework and the techniques developed in this work would also be an interesting direction of future work.

\bibliographystyle{IEEEtran}
\bibliography{IEEEabrv,CentralBib,glrt}
\begin{IEEEbiography}[{\includegraphics[width=1in,height=1.25in,clip,keepaspectratio]{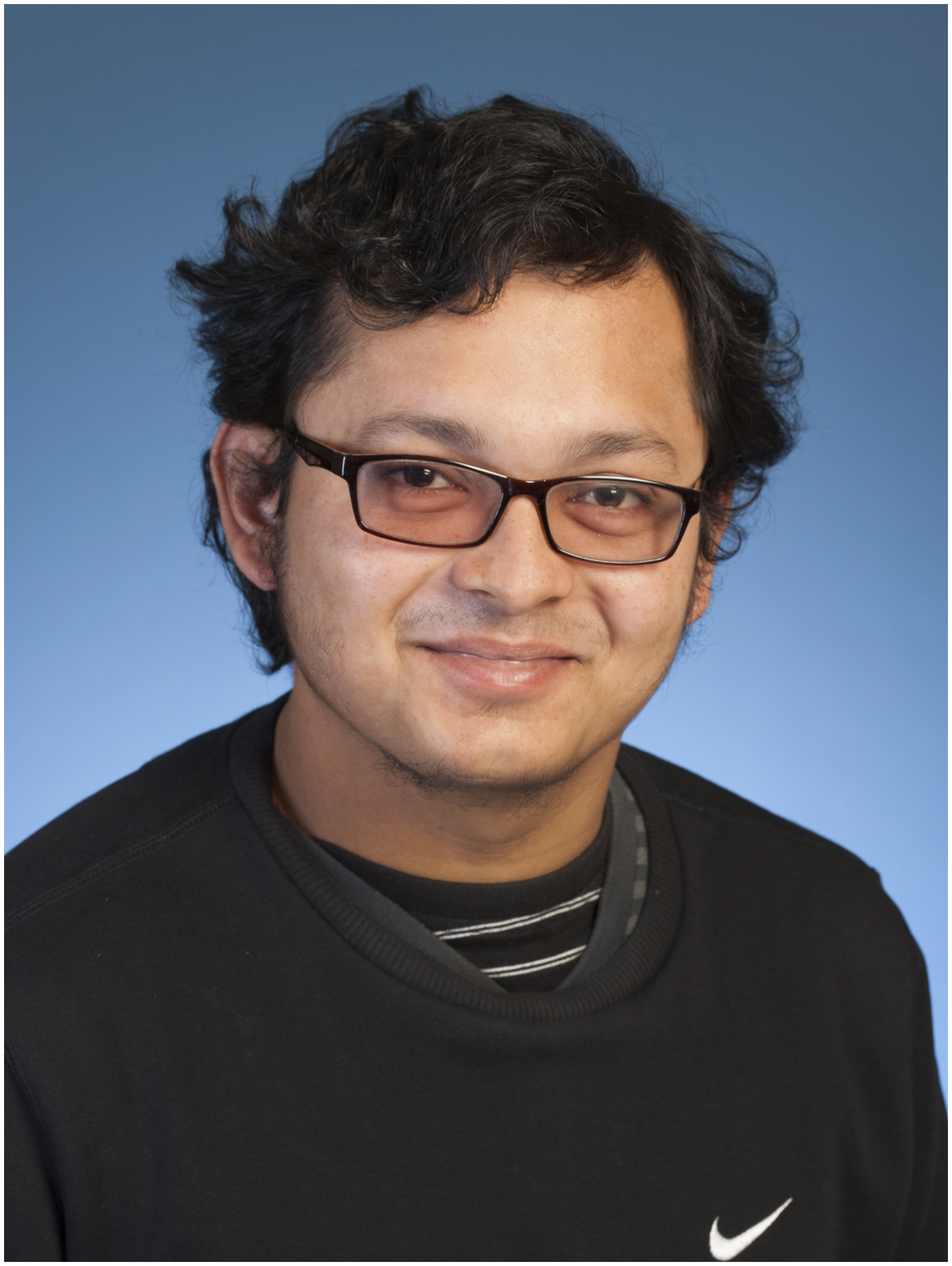}}]{Anit Kumar Sahu} (S'13) received
	a B.Tech. in Electronics and Electrical
	Communication Engineering and M.Tech. in Telecommunication Systems Engineering from the Indian Institute
	of Technology, Kharagpur, India, in May 2013. Since Fall 2013, he has been working towards his Ph.D. in Electrical and Computer Engineering at Carnegie Mellon University, Pittsburgh, PA. His research interests include distributed inference in large-scale stochastic systems, statistical machine learning, and information theory.
\end{IEEEbiography}

\begin{IEEEbiography}[{\includegraphics[width=1in,height=1.25in,clip,keepaspectratio]{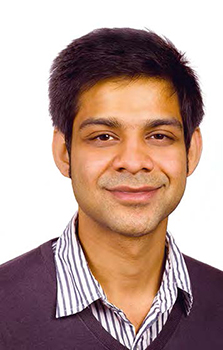}}]{Soummya Kar} (S'05, M'10) received a B.Tech. in electronics and electrical communication 
	engineering from the Indian Institute of Technology, Kharagpur, India, in May 2005 and 
	a Ph.D. in electrical and computer engineering from Carnegie Mellon University, 
	Pittsburgh, PA, in 2010. From June 2010 to May 2011, he was with the Electrical Engineering Department, Princeton University, Princeton, NJ, USA, as a Postdoctoral Research Associate. He is currently an Associate Professor of Electrical and Computer Engineering at Carnegie Mellon University, Pittsburgh, PA, USA. His research interests include decision-making in large-scale networked  systems, stochastic systems, multi-agent systems and data science, with applications to cyber-physical systems and smart energy systems. He has published extensively in these topics with more than 140 articles in journals and conference proceedings and holds multiple patents.  Recent recognition of his work includes the 2016 O. Hugo Schuck Best Paper Award from the American Automatic Control Council, the 2016 Dean's Early Career Fellowship from CIT, Carnegie Mellon, and the 2011 A.G. Milnes Award for best PhD thesis in Electrical and Computer Engineering, Carnegie Mellon University.
\end{IEEEbiography}
\begin{IEEEbiography}[{\includegraphics[width=1in,height=1.25in,clip,keepaspectratio]{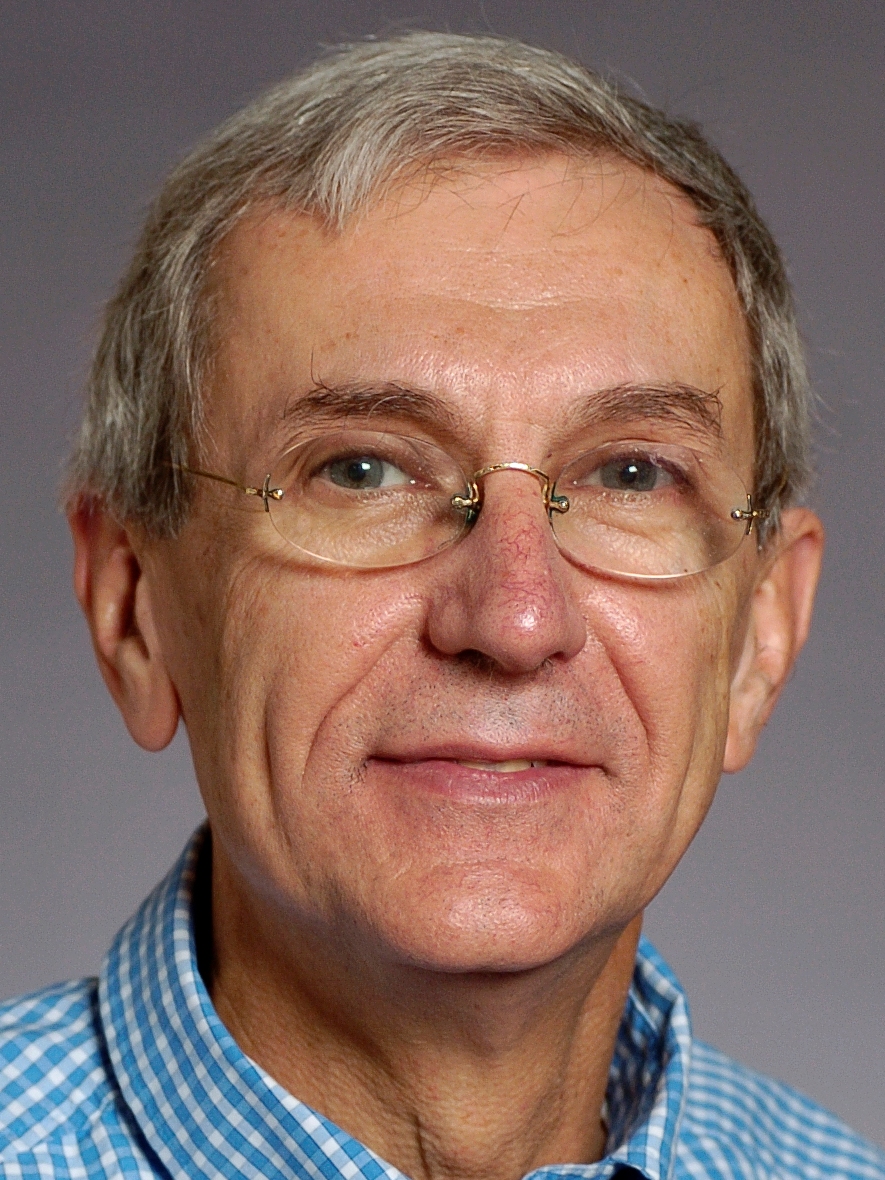}}]{Jos\'e M.~F.~Moura}(S'71--M'75--SM'90--F'94) received the engenheiro electrot\'{e}cnico degree from Instituto Superior T\'ecnico (IST), Lisbon, Portugal, and the M.Sc., E.E., and D.Sc.~degrees in EECS from the Massachusetts Institue of Technology (MIT), Cambridge, MA.
	
	He is the Philip L. and Marsha Dowd University Professor at Carnegie Mellon University (CMU). He was on the faculty at IST and has held visiting faculty appointments at MIT and New York University (NYU). He founded and directs a large education and research program between CMU and Portugal, www.icti.cmu.edu.
	
	His research interests are on  data science, graph signal processing, and statistical and algebraic signal and image processing. He has published over 550 papers and holds thirteen patents issued by the US Patent Office. The technology of two of his patents (co-inventor A. Kav\v{c}i\'c) are in about three billion disk drives read channel chips of 60~\% of all computers sold in the last 13 years worldwide and were, in 2016, the subject of the largest university verdict/settlement in the information technologies area.
	
	Dr. Moura is the IEEE Technical Activities Vice-President (2016) and member of the IEEE Board of Directors. He served in several other capacities including IEEE Division IX Director, member of several IEEE Boards, President of the IEEE Signal Processing Society(SPS), Editor in Chief for the IEEE Transactions in Signal Processing, interim Editor in Chief for the IEEE Signal Processing Letters.
	
	Dr. Moura has received several awards, including  the Technical Achievement Award and the Society Award from the IEEE Signal Processing. In 2016, he received the CMU College of Engineering Distinguished Professor of Engineering Award. He is a Fellow of the IEEE, a Fellow of the American Association for the Advancement of Science (AAAS), a corresponding member of the Academy of Sciences of Portugal, Fellow of the US National Academy of Inventors, and a member of the US National Academy of Engineering.
\end{IEEEbiography}

\begin{IEEEbiography}[{\includegraphics[width=1in,height=1.25in,clip,keepaspectratio]{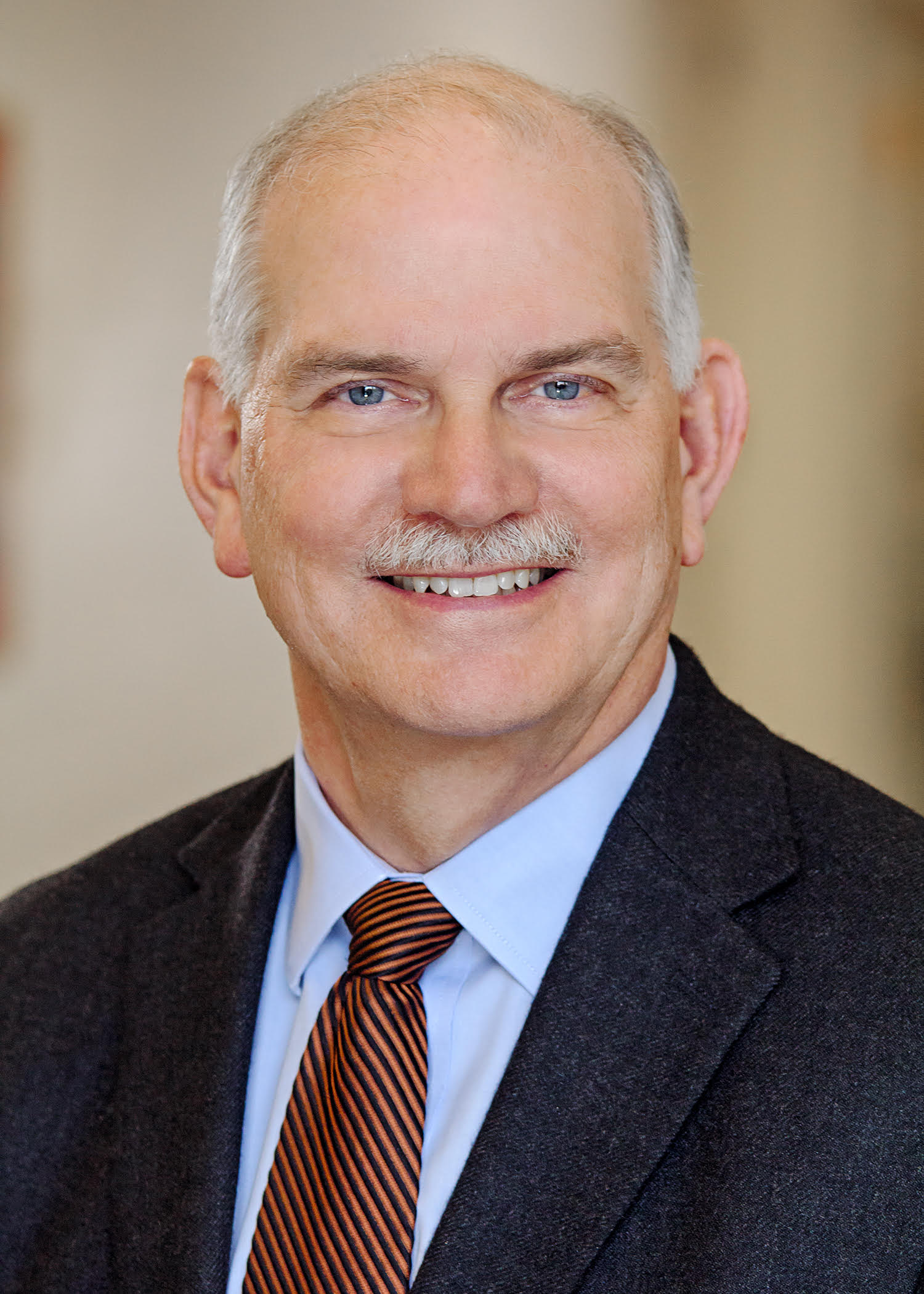}}]{H. Vincent Poor} (S'72, M'77, SM'82, F'87) received the Ph.D. degree in EECS from Princeton University in 1977.  From 1977 until 1990, he was on the faculty of the University of Illinois at Urbana-Champaign. Since 1990 he has been on the faculty at Princeton, where he is the Michael Henry Strater University Professor of Electrical Engineering. During 2006 to 2016, he served as Dean of Princeton’s School of Engineering and Applied Science. Dr. Poor’s research interests are in the areas of statistical signal processing, stochastic analysis and information theory, and their applications in wireless networks and related fields. Among his publications in these areas is the recent book \emph{Mechanisms and Games for Dynamic Spectrum Allocation} (Cambridge University Press, 2014).\\

\noindent Dr. Poor is a member of the National Academy of Engineering and the National Academy of Sciences, and a foreign member of the Royal Society. He is also a Fellow of the American Academy of Arts and Sciences and the National Academy of Inventors, and of other national and international academies. He received the Technical Achievement and Society Awards of the IEEE Signal Processing Society in 2007 and 2011, respectively. Recent recognition of his work includes the 2014 URSI Booker Gold Medal, the 2015 EURASIP Athanasios Papoulis Award, the 2016 John Fritz Medal, and honorary doctorates from Aalborg University, Aalto University, HKUST and the University of Edinburgh.
\end{IEEEbiography}

\end{document}